# ANALYSIS AND SIMULATIONS ON A MODEL FOR THE EVOLUTION OF TUMORS UNDER THE INFLUENCE OF NUTRIENT AND DRUG APPLICATION

KONSTANTINA TRIVISA AND FRANZISKA WEBER

ABSTRACT. We investigate the evolution of tumor growth relying on a nonlinear model of partial differential equations which incorporates mechanical laws for tissue compression combined with rules for nutrients availability and drug application. Rigorous analysis and simulations are presented which show the role of nutrient and drug application in the progression of tumors. We construct an explicit convergent numerical scheme to approximate solutions of the nonlinear system of partial differential equations. Extensive numerical tests show that solutions exhibit a necrotic core when the nutrient level falls below a critical level in accordance with medical observations. The same numerical experiment is performed in the case of drug application for the purpose of comparison. Depending on the balance between nutrient and drug both shrinkage and growth of tumors can occur. The role of inhomogeneous boundary conditions, vascularization and anisotropies in the development of tumor shape irregularities are discussed.

1. INTRODUCTION

1.1. **Motivation.** Tumor cells can show distinct morphological and phenotypic profiles, including cellular morphology, gene expression, metabolism, motility, proliferation, and metastatic potential. The heterogeneity of cancer cells introduces significant challenges in designing effective treatment strategies. Scientific research aiming at understanding and characterizing heterogeneity can allow for a better understanding of the causes and progression of disease. In turn, this has the potential to guide the creation of more refined treatment strategies that incorporate knowledge of heterogeneity to yield higher efficacy. In recent years, the investigation of the effect of drug application in cancer progression has been the subject of intense scientific effort. Both vascularization and anisotropies affect the evolution of tumors, the shape of the surface of the tumor regime and the development of irregularities. Mathematical modeling, the construction of explicit examples and the use of simulations can establish critical conditions for the growth or the shrinkage of tumors, the development of irregularities for different types of cancers enhancing our understanding of tumor developement and heterogeneity.

1.2. **Governing equations.** Over the past years research activity on the mathematical modeling and simulations on tumor growth models has increased dramatically. A variety of modeling strategies have been developed, each focusing on one or more aspects of cancer for medical prediction. Among the variety of models now available one needs to mention as a starting point the class, introduced by

*Date*: February 21, 2016.
2010 *Mathematics Subject Classification*. Primary: 35Q30, 76N10; Secondary: 92C17.
*Key words and phrases*. Tumor growth models, cancer progression, finite difference schemes, multi-phase flow, existence of solutions.
The work of K.T. was supported in part by the National Science Foundation under the grant DMS-1211519.

Greenspan [9], which considers that cancerous cells multiplication is due to nutrients (glucosis, oxygen) brought by blood vessels. The first stage, where growth is limited by nutrients, lasts until the tumor reaches a certain size; subsequently, lack of food leads to cell necrosis which triggers neovasculatures development that supply the tumor with enough nourishment. This process has motivated a new generation of models where growth is effected both by the nutrient supply and the competition for space, enriching the modeling effort with mechanical considerations, viewing tissues as multiphase fluids. This is the point of view adopted in our investigation.

1.2.1. *Transport equations for the evolution of the cell densities.* All the cells are assumed to follow the general continuity equation:

$$\partial_t n - \operatorname{div}(n\boldsymbol{u}) = n\boldsymbol{\Phi}(p,c,q), \quad x \in \Omega,\ t \geq 0 \tag{1.1}$$

where $n$ represents the number density of tumor cells, $\boldsymbol{u}$ the velocity field, $c$ is the concentration of the nutrient (oxygen) and $q$ the density of the drug. In the present context, $p$ denotes the pressure of the *tumor* and $\Omega$ a bounded domain in $\mathbb{R}^d$, $d=2,3$. The function $\boldsymbol{\Phi}$ accounts for the effect of pressure, nutrient and drug to the evolution of cancerous cells and has the general form

$$\boldsymbol{\Phi}(p,c,q) = g_1(c,q)\boldsymbol{G}(p) - g_2(c,q), \tag{1.2}$$

where $g_i$, $i=1,2$ are bounded and nonnegative. For $\boldsymbol{G}$ we assume that it is of the form

$$\boldsymbol{G} \in C^1(\mathbb{R}), \quad \boldsymbol{G}'(\cdot) \leq 0, \quad \boldsymbol{G}(p_M) = 0 \quad \text{for some } P_M > 0, \tag{1.3}$$

where $P_M$ is the so called *homeostatic pressure*, the critical threshold at which the cell division is stopped by contact inhibition. It is related to the compression a cell can experience [3]. The factor $g_1(c,q)$ is nondecreasing in $c$ and nonincreasing in $q$ and accounts for the growth/decline of the cell culture in relation to the nutrient and drug concentration. The term $g_2(c,q)$ accounts for the decrease of the cell density when the nutrient concentration falls below a critical level to sustain cell life (cells starve) and also models the effect of the drug application on the tumor growth. It can for example be of the form $g_2(c,q) = g_{2,1}(c) + g_{2,2}(q)$. The pressure law is given by

$$p(n) = n^\gamma, \tag{1.4}$$

where $\gamma \geq 2$. Here, and in what follows, for simplicity we let

$$\boldsymbol{G}(p) = \alpha - \beta \widetilde{g}_0(p), \tag{1.5}$$

for some $\alpha, \beta > 0$, where $\widetilde{g}_0(p) = g_0(n)$ is such that $ng_0(n)$ is nondecreasing. This includes for example functions of the form $\widetilde{g}_0(p) = p^\theta$ for $\theta > 0$.

1.2.2. *The tumor tissue as a porous medium.* The continuous motion of cells within the tumor region typically due to proliferation is represented by the velocity field $\boldsymbol{u} := \nabla W$ given by an alternative to Darcy's equation known as *Brinkman's equation*

$$p = W - \mu \Delta W \tag{1.6}$$

where $\mu$ is a positive constant describing the viscous like properties of tumor cells and $p$ is the pressure given by (1.4).

Relation (1.6) extents the usual Darcy's law, by taking into account the dissipative force density, which results from the internal cell friction due to cell volume changes.

1.2.3. *A linear diffusion equation for the evolution of nutrient.* Tumor cells consume nutrients (oxygen). In contrast to the equations of cell densities, the equations of the oxygen molecules in the tumor include diffusion terms in the following form:
$$\partial_t c - \nabla \cdot (\nu_c \nabla c) = c\Psi_c(n,c) + r_c(c_{\text{supp}} - c).$$
Assuming that $\nu_c$ is constant this equation (cf. Friedman [8]), the equation becomes
$$\begin{cases} \partial_t c - \nu_c \Delta c = c\Psi_c(n,c) + r_c(c_{\text{supp}} - c). \\ c(t,x) = c_b(t,x) \geq 0 \text{ for } x \in \partial\Omega. \end{cases} \quad (1.7)$$

The function $\Psi_c$ is nonpositive and accounts for the consumption of the nutrient by the cell culture. According to (cf. Ward and King [18, 17]) the nutrient is consumed at a rate proportional to the rate of *cell mitosis*, which is accounted by the first term on the right-hand side. $r_c > 0$ is the rate at which the nutrient is supplied to the tumor region (term with $c_{\text{supp}}$) and the consumption by he healthy cells (term with $c$).

1.2.4. *A linear diffusion equation for the evolution of drug.* The evolution of the drug concentration in the tumor is given by a diffusion equation of the form
$$\partial_t q - \nabla \cdot (\nu_q \nabla q) = q\Psi_q(n,q) + r_q(q_{\text{supp}} - q).$$
Assuming that $\nu_q$ is constant this equation becomes
$$\begin{cases} \partial_t q - \nu_q \Delta q = q\Psi_q(n,q) + r_q(q_{\text{supp}} - q), \\ q(t,x) = q_b(t,x) \geq 0, \text{ for } x \in \partial\Omega. \end{cases} \quad (1.8)$$

with $\Psi_q(\cdot)$ a smooth nonpositive function. This equation describes the diffusion of the drug within the tumor region. The first term of the right-hand side of (1.8) represents the drug consumption and can be viewed as a measure of the drug effectiveness. The second term on the right represent the rate by which the drug is applied to the tumor region. The resulting model, governed by the transport equation (1.1) for the population density of cells, the elliptic equation (1.6) for the velocity field and a state equation for the pressure law (1.4), now reads

$$\partial_t n - \text{div}(n \nabla W) = n\Phi(p,c,q), \quad x \in \Omega, \ t \geq 0 \quad (1.9a)$$
$$-\mu \Delta W + W = n^\gamma \quad (1.9b)$$
$$\partial_t c - \nu_c \Delta c = c\Psi_c(n,c) + r_c(c_{\text{supp}} - q) \quad (1.9c)$$
$$\partial_t q - \nu_q \Delta q = q\Psi_q(n,q) + r_q(q_{\text{supp}} - q). \quad (1.9d)$$

We complete the system (1.9) with a family of initial data $n_0$ satisfying (for some constant $C$)
$$\begin{cases} n_0 \geq 0, \quad p(n_0) \leq P_M, \quad \|n_0\|_{L^1(\mathbb{R}^d)} \leq C \\ 0 \leq c_0 \leq c_\infty, \quad 0 \leq q_0 \leq q_\infty. \end{cases} \quad (1.10)$$
for some constants $0 \leq c_\infty, q_\infty < \infty$.

The objective of this work is to design an efficient numerical scheme for the approximation of the solution to the nonlinear system (1.9) and to establish that this scheme converges when the mesh is refined; yielding at the same time the global existence of weak solutions to the nonlinear model for tumor growth (1.9). The main ingredients of our approach and contribution to the existing theory on Hele-Shaw-type systems for tumor growth include:

- The design of an efficient numerical scheme for the numerical approximation of the nonlinear system (1.9a)-(1.9d) with the aid of a finite difference scheme.

- The proof of the convergence of the numerical scheme, which is achieved by establishing the strong convergence of the cell densities. This property is obtained as a consequence of the weak continuity of the *effective viscous pressure* (cf. Section 4).
- The construction of numerical experiments that establish various cancer phenomena confirmed by clinical observations providing criteria for the development of certain types of tumor heterogeneities.

For relevant results on the analysis and the numerical approximation of a two-phase flow model in porous media we refer the reader to [5]. Related work on the mathematical analysis of mechanical models of Hele-Shaw-type have been presented by Perthame *et al.* [14, 15]. In [16], Trivisa and Weber presented a convergent explicit finite difference scheme for the numerical approximation of a Hele-Shaw-type system for the evolution of cancerous cells and presents numerical observations in two space dimensions. The work [16] is according to our knowledge the first article that presents rigorous analytical results on the global existence of general weak solutions to Hele-Shaw-type systems. The present article extends the analysis in [16] significantly by investigating the delicate interplay of nutrient and drug application in the treatment of cancer. We refer the reader to [7, 8, 4] where a class of relevant tumor growth models with nutrient and drug application are presented and to the manuscript [1] which provides an overview of mathematical methods and tools for modeling cancer phenomena.

1.3. **Outline.** The paper is organized as follows: Section 1 presents the motivation, modeling and introduces the necessary preliminary material. Section 2 provides a weak formulation of the problem and states the main results. Section 3 is devoted to the global existence of solutions via a vanishing viscosity approximation. In Section 4 we present an efficient finite difference scheme for the approximation of the weak solution to system (1.9) on rectangular domains and prove its convergence. We conclude by presenting extensive numerical tests in Section 5.

## 2. Weak formulation and main results

We start by defining a notion of a weak solution to system (1.9):

**Definition 2.1.** Let $\Omega$ a bounded domain in $\mathbb{R}^d$, $d = 2, 3$, with smooth boundary $\partial\Omega$ and $T > 0$ a finite time horizon. We say that $(n, W, p, c, q)$ is a weak solution of problem (1.1)-(1.6) supplemented with initial data $(n_0, W_0, p_0, c_0, q_0)$ satisfying (1.10) provided that the following hold:

- $(n, W, p, c, q) \geq 0$ represents a weak solution of (1.1)–(1.6) on $(0, T) \times \Omega$, i.e., for any test function $\varphi \in C_c^\infty([0, T] \times \mathbb{R}^d), T > 0$, the following integral relations hold

$$\int_{\mathbb{R}^d} n\varphi(\tau, \cdot)\, dx - \int_{\mathbb{R}^d} n_0 \varphi(0, \cdot) dx =$$
$$\int_0^\tau \int_{\mathbb{R}^d} (n\partial_t \varphi - n\nabla W \cdot \nabla\varphi + n\boldsymbol{\Phi}(p, c, q)\varphi(t, \cdot))\, dxdt, \quad (2.1)$$

and

$$n \in L^q((0, T) \times \Omega), \text{ for all } q \geq 1.$$

We remark that in the weak formulation, it is convenient that the equations (1.1) hold on the whole space $\mathbb{R}^d$ provided that the densities $n$ are extended to be zero outside the tumor domain.

• Brinkman's equation (1.6) holds in the sense of distributions, i.e., for any test function $\varphi \in C_c^\infty(\mathbb{R}^d)$, the following integral relation holds for a.e. $t \in [0,T]$,

$$\int_\Omega n^\gamma \varphi \, dx = \int_\Omega \left(\mu \nabla W \cdot \nabla \varphi + W\varphi\right) dx. \qquad (2.2)$$

and $p = n^\gamma$ almost everywhere. All quantities in (2.2) are required to be integrable, in particular, $W \in L^\infty([0,T]; H^2(\Omega))$.

• The diffusion equation for the evolution of the density of the nutrient (1.7) holds in the sense of distribution, i.e., for any test function $\varphi \in C_c^\infty([0,T) \times \mathbb{R}^d), T > 0$ with compact support in $[0,T) \times \Omega$, we have

$$\int_{\mathbb{R}^d} c\varphi(\tau, \cdot) \, dx - \int_{\mathbb{R}^d} c_0 \varphi(0,\cdot) dx = \int_0^\tau \int_{\mathbb{R}^d} c\partial_t \varphi \, dx dt -$$
$$\int_0^\tau \int_{\mathbb{R}^d} \nu_c \nabla_x c \cdot \nabla_x \varphi dx dt + \int_0^\tau \int_{\mathbb{R}^d} (c\boldsymbol{\Psi}_c(n,c) + r_c(c_{\text{supp}} - c))\varphi dx dt,$$

and $c \in L^\infty((0,T) \times \Omega) \cap L^2([0,T]; H^1(\Omega))$.

• $q \geq 0$ is a weak solution of (1.8), i.e., for any test function $\varphi \in C_c^\infty([0,T) \times \mathbb{R}^d), T > 0$ with compact support in $[0,T) \times \Omega$, the following integral relations hold

$$\int_{\mathbb{R}^d} q\varphi(\tau,\cdot) \, dx - \int_{\mathbb{R}^d} q_0 \varphi(0,\cdot) dx = \int_0^\tau \int_{\mathbb{R}^d} q\partial_t \varphi dx dt -$$
$$\int_0^\tau \int_{\mathbb{R}^d} \nu_q \nabla_x q \cdot \nabla_x \varphi dx dt + \int_0^\tau \int_{\mathbb{R}^d} (q\boldsymbol{\Psi}_q(n,q) + r_q(q_{\text{supp}} - q))\varphi dx dt,$$

and $q \in L^\infty((0,T) \times \Omega) \cap L^2([0,T]; H^1(\Omega))$.

The main result of the article now follows.

**Theorem 2.2.** *Let $\Omega \subset \mathbb{R}^d$ be a bounded domain with smooth boundary $\partial \Omega$, $0 < T < \infty$. Assume that the initial data $n_0 \in L^\infty(\Omega)$ with $0 \leq n_0 \leq n_\infty := P_M^{1/\gamma}$ and that $\boldsymbol{\Phi}(\cdot)$ is of the form (1.2), (1.3) and (1.5). Then the problem (1.1)–(1.6) admits a weak solution in the sense specified in Definition 2.1.*

In Section 3, we will obtain such a solution as the limit of a vanishing viscosity approximation $(n_\varepsilon, W_\varepsilon, p_\varepsilon, c_\varepsilon, q_\varepsilon)$ of (3.1) to (1.9) as $\varepsilon \to 0$ on a domain with smooth boundary $\partial \Omega$ and in Section 4 as the limit of the sequence of approximations $(n_h, W_h, p_h, c_h, q_h)$ computed by the numerical scheme (4.1) – (4.3) as $h \to 0$ on a rectangular domain $\Omega$.

## 3. Global existence via vanishing viscosity

In this section we prove Theorem 2.2 by constructing an approximating scheme which relies on the addition of artificial viscosity in the cell density evolution equation

$$\partial_t n_\varepsilon - \operatorname{div}(n_\varepsilon \nabla W_\varepsilon) = n_\varepsilon \boldsymbol{\Phi}(p_\varepsilon, c_\varepsilon, q_\varepsilon) + \varepsilon \Delta n_\varepsilon, \quad x \in \Omega, \ t \geq 0 \qquad (3.1\text{a})$$
$$\mu \Delta W_\varepsilon - W_\varepsilon = n_\varepsilon^\gamma, \qquad (3.1\text{b})$$
$$\partial_t c_\varepsilon - \nu_c \Delta c_\varepsilon = c_\varepsilon \boldsymbol{\Psi}_c(n_\varepsilon, c_\varepsilon) + r_c(c_{\text{supp}} - c_\varepsilon) \qquad (3.1\text{c})$$
$$\partial_t q_\varepsilon - \nu_q \Delta q_\varepsilon = q_\varepsilon \Psi_q(n_\varepsilon, q_\varepsilon) + r_q(q_{\text{supp}} - q_\varepsilon), \qquad (3.1\text{d})$$
$$n_\varepsilon(0,\cdot) = n_0^\varepsilon, \quad c_\varepsilon(0,\cdot) = c_0^\varepsilon, \quad q_\varepsilon(0,\cdot) = q_0^\varepsilon, \qquad (3.1\text{e})$$

where $n_0^\varepsilon$, $c_0^\varepsilon$ and $q_0^\varepsilon$ are smoothnened versions of $n_0$, $c_0$, $q_0$ respectively, that is $n_0^\varepsilon = n_0 * \varphi_\varepsilon$ for a smooth function $\varphi_\varepsilon$ with compact support. $\Omega$ is a bounded domain in $\mathbb{R}^d$ with smooth boundary or alternatively the $d$-dimensional torus $\mathbb{T}^d$. We assume homogeneous Neumann boundary conditions for $n_\varepsilon$ and $W_\varepsilon$ (if the domain is a torus $\mathbb{T}^d$ we can also use periodic boundary conditions). For $c$ and $q$ we assume that they satisfy a Dirichlet boundary condition of the form

$$c(t,x) = c_b(t,x), \quad q(t,x) = q_b(t,x), \quad x \in \partial\Omega, \, t \in [0,T],$$

for smooth $c_b$ and $q_b$ with $c_b(0,x) = c_0(x)$ and $q_b(0,x) = q_0(x)$ for $x \in \partial\Omega$. We assume that $c_b$ and $q_b$ are such that there exists smooth extensions $\widetilde{c}_b, \widetilde{q}_b \in C^2([0,T] \times \Omega)$ of them onto the whole domain $[0,T] \times \overline{\Omega}$. In addition, we assume that for constants $0 < c_\infty, q_\infty < \infty$,

$$0 \leq c_b(t,x), c_0(x), c_{\text{supp}} \leq c_\infty, \quad 0 \leq q_b(t,x), q_0(x), c_{\text{supp}} \leq q_\infty.$$

**Theorem 3.1.** *For every $\varepsilon > 0$, the parabolic-elliptic system* (3.1) *admits a unique smooth solution* $(n_\varepsilon, W_\varepsilon, p_\varepsilon, c_\varepsilon, q_\varepsilon)$.

*Proof.* The proof of this result relies on classical arguments (cf. Ladyzhenskaya [10]). For details we refer the reader to Lunardi [12, Theorem 5.1.2] in the context of a related parabolic partial differential equation. □

The remaining part of this section aims to establish the necessary compactness of the approximate sequence of solutions $(n_\varepsilon, W_\varepsilon, p_\varepsilon, c_\varepsilon, q_\varepsilon)$.

3.1. **A priori estimates.** We start by proving that $n_\varepsilon$ are uniformly bounded independent of $\varepsilon > 0$ and nonnegative:

**Lemma 3.2.** *For any $t > 0$, the functions $n_\varepsilon(t,\cdot)$ are uniformly (in $\varepsilon > 0$) bounded and nonnegative if $0 \leq n_\varepsilon(0,\cdot) \leq n_\infty := P_M^{1/\gamma} < \infty$ uniformly in $\varepsilon > 0$, specifically*

$$0 \leq \min_{(t,x)} n_\varepsilon(t,x) \leq \max_{(t,x)} n_\varepsilon(t,x) \leq n_\infty.$$

*Proof.* The proof of this lemma is similar to the proof of Lemma 3.2 in [16]. All that needs to be done is replacing the function $\boldsymbol{G}$ on the right hand side of the equation by $\boldsymbol{\Phi}$ and checking that it satisfies the right growth conditions. □

As a next step, we prove positivity and uniform boundedness of the nutrient sequence $\{c_\varepsilon\}_{\varepsilon \geq 0}$.

**Lemma 3.3.** *For any $t > 0$, the functions $c_\varepsilon(t,\cdot)$ are uniformly (in $\varepsilon > 0$) bounded and nonnegative if $0 \leq c_{supp}, c_\varepsilon(0,\cdot), c_b(t,x) \leq c_\infty < \infty$ for all $\varepsilon > 0$, specifically*

$$0 \leq \min_{(t,x)} c_\varepsilon(t,x) \leq \max_{(t,x)} c_\varepsilon(t,x) \leq c_\infty.$$

*Proof.* For the ease of notation, we omit writing the subscript $\varepsilon$. To prove the nonnegativity, let us assume that $(t_0, x_0)$ is a point where $c(t_0, x_0) = 0$ for the first time (i.e. $c(t,x) > 0$ for any $0 \leq t < t_0$ and any $x \in \Omega$). Then $c(t_0, x) \geq 0$ for all $x$ in some small neigborhood of $x_0$ and hence $\Delta c(t_0, x_0) \geq 0$. Thus, bringing the Laplace term to the right hand side of equation (3.1c), we have

$$\partial_t c(t_0, x_0) = \nu_c \Delta c(t_0, x_0) + c(t_0, x_0) \boldsymbol{\Psi}_c(n, c(t_0, x_0)) + r_c(c_{\text{supp}} - c(t_0, x_0))$$
$$= \nu_c \Delta c(t_0, x_0) + r_c c_{\text{supp}} \geq r_c c_{\text{supp}} > 0,$$

where we used that by assumption $c(t_0, x_0) = 0$ for the second equality and that $\Delta c(t_0, x_0) \geq 0$ for the inequality. Thus at a point where $c$ becomes zero, the time derivative of $c$ is positive and it will thus stay nonnegative. To show the uniform boundedness, let us denote $\widehat{c} := c_\infty - c$. Then, $\widehat{c}$ satisfies the equation

$$\partial_t \widehat{c} = \nu_c \Delta \widehat{c} - (c_\infty - \widehat{c}) \boldsymbol{\Psi}_c(n, c_\infty - \widehat{c}) + r_c(c_\infty - c_{\text{supp}} - \widehat{c}).$$

We show that $\widehat{c}$ remains nonnegative (it is initially nonnegative by the assumptions on the initial data $c_0$ and on the boundary by assumptions on $c_b$). We therefore again assume that $(t_0, x_0)$ is a point where $\widehat{c}(t_0, x_0) = 0$ for the first time (i.e. $\widehat{c}(t, x) > 0$ for any $0 \leq t < t_0$ and any $x \in \Omega$). Then $\widehat{c}(t_0, x) \geq 0$ for all $x$ in some small neigborhood of $x_0$ and hence $\Delta \widehat{c}(t_0, x_0) \geq 0$. We use this fact and that $\widehat{c}(t_0, x_0) = 0$ in the evolution equation for $\widehat{c}$:

$$\partial_t \widehat{c}(t_0, x_0) \geq -c_\infty \mathbf{\Psi}_c(n, c_\infty) + r_c(c_\infty - c_{\text{supp}}).$$

By the assumptions on $\mathbf{\Psi}_c$, $\mathbf{\Psi}_c(n, c_\infty) \leq 0$, and hence the right hand side is non-negative. Therefore, $\widehat{c}$ remains nonnegative which implies the boundedness of $c$. $\square$

Positivity and uniform boundedness of the sequence of the drug functions $\{q_\varepsilon\}_{\varepsilon \geq 0}$ is proved in the same way. Next, we prove a regularity estimate for the nutrient $c_\varepsilon$.

**Lemma 3.4.** *We have that, uniformly for all $\varepsilon > 0$,*

$$c_\varepsilon \in L^2([0,T]; H^1(\Omega)), \partial_t c_\varepsilon \in L^2([0,T]; H^{-1}(\Omega)).$$

*In particular,*

$$\int_0^T \int_\Omega |\nabla c_\varepsilon|^2 \, dx dt \leq C,$$

*where $C$ is a constant independent of $\varepsilon > 0$.*

*Proof.* For the ease of notation, we omit writing the subscript $\varepsilon$ of $c_\varepsilon$. Let us define the function $\widehat{c} := c - \widetilde{c}_b$, where $\widetilde{c}_b$ is the smooth extension of $c_b$ to the whole domain. It satisfies the equation

$$\partial_t \widehat{c} - \nu_c \Delta \widehat{c} = (\widehat{c} + \widetilde{c}_b) \mathbf{\Psi}_c(n, \widehat{c} + \widetilde{c}_b) + r_c(c_{\text{supp}} - \widehat{c} - \widetilde{c}_b) - \partial_t \widetilde{c}_b + \Delta \widetilde{c}_b, \quad x \in \Omega, t \in [0, T]$$

$$\widehat{c}(x, t) = 0, \quad x \in \partial \Omega, t \in [0, T].$$

Now we multiply the evolution equation for $\widehat{c}$ by $\widehat{c}$ and integrate over the spatial domain. After integration by parts and using the homogeneous Dirichlet boundary conditions for the boundary integral contributions, we obtain,

$$\frac{1}{2} \frac{d}{dt} \int_\Omega \widehat{c}^2 \, dx + \nu_c \int_\Omega |\nabla \widehat{c}|^2 \, dx = \int_\Omega [\widehat{c}(\widehat{c} + \widetilde{c}_b) \mathbf{\Psi}_c(n, \widehat{c} + \widetilde{c}_b) + r_c \widehat{c}(c_{\text{supp}} - \widehat{c} - \widetilde{c}_b)] \, dx$$

$$+ \int_\Omega \widehat{c} (\Delta \widetilde{c}_b - \partial_t \widetilde{c}_b) \, dx$$

Thanks to the $L^\infty$-bounds on $n_\varepsilon$ and $c_\varepsilon$ from Lemmas 3.2 and 3.3 and the smoothness of $\widetilde{c}_b$, we obtain that the right hand side is bounded. Thus integrating in time, we obtain that $\widehat{c} \in L^2([0,T]; H^1(\Omega))$. Since $c$ differs from $\widehat{c}$ only by the smooth function $\widetilde{c}_b$, we obtain that the same holds for $c$. The $L^2([0,T]; H^{-1}(\Omega))$-bound on $\partial_t c$ then follows using this, from the equation it satisfies. $\square$

*Remark* 3.5. It is possible to obtain higher order interior regularity estimates for $c_\varepsilon$, specifically, one can show $c_\varepsilon \in L^\infty([0,T]; H^1(V)) \cap L^2([0,T]; H^2(V))$, $\partial_t c_\varepsilon \in L^2([0,T] \times V)$ for any compact subset $V \subset\subset \Omega$. However, we will not need this for the proof of convergence of the approximating sequence and therefore omit the proof of this fact here.

*Remark* 3.6. Using exactly the same arguments, it is shown that $q_\varepsilon \in L^2([0,T]; H^1(\Omega))$ and $\partial_t q_\varepsilon \in L^2([0,T]; H^{-1}(\Omega))$.

Next we recall Lemma 3.3 from [16] which gives us higher order regularity estimates for $W_\varepsilon$.

**Lemma 3.7.** *We have that*

$$W_\varepsilon \in L^\infty([0,T]; H^2(\Omega)), \quad W_\varepsilon \in L^\infty([0,T]; W^{2,q}(\Omega')),$$

*for any $q \in [1, \infty)$, all compact subsets $\Omega' \subset\subset \Omega$, uniformly in $\varepsilon > 0$ and*

$$W_\varepsilon, \Delta W_\varepsilon \in L^\infty((0,T) \times \Omega),$$

*uniformly in $\varepsilon > 0$ as well.*

*Proof.* This lemma was proved as Lemma 3.3 in [16]. □

### 3.2. Entropy inequalities for $n_\varepsilon$.

To prove strong convergence of the approximating sequence $\{(n_\varepsilon, W_\varepsilon, p_\varepsilon, c_\varepsilon, q_\varepsilon)\}_{\varepsilon > 0}$, it is useful to derive entropy inequalities for $n_\varepsilon$. To this end, we recall the following lemma:

**Lemma 3.8.** *Let $f : \mathbb{R} \to \mathbb{R}$ be a smooth convex function and denote $f_\varepsilon := f(n_\varepsilon)$. Then $f_\varepsilon$ satisfies the following identity*

$$\partial_t f_\varepsilon - \mathrm{div}(f_\varepsilon \nabla W_\varepsilon) - \varepsilon \Delta f(n_\varepsilon)$$
$$= (f'(n_\varepsilon) n_\varepsilon - f_\varepsilon) \Delta W_\varepsilon + f'(n_\varepsilon) n_\varepsilon \mathbf{\Phi}(p_\varepsilon, c_\varepsilon, q_\varepsilon) - \varepsilon f''(n_\varepsilon) |\nabla n_\varepsilon|^2 \quad (3.2)$$

*where*

$$\varepsilon \int_0^T \int_\Omega f''(n_\varepsilon) |\nabla n_\varepsilon|^2 \, dx dt \leq C, \quad (3.3)$$

*with $C > 0$ a constant independent of $\varepsilon > 0$. In particular, this implies that $\partial_t f_\varepsilon = g_\varepsilon + k_\varepsilon$ with $g_\varepsilon \in L^1([0,T] \times \Omega)$ and $k_\varepsilon \in L^1([0,T]; W^{-1,2}(\Omega))$.*

*Proof.* The identity (3.2) follows by multiplying the evolution equation for $n_\varepsilon$, (3.1a), by $f'(n_\varepsilon)$, applying chain rule, integrating the inequality in space and time, and following similar line of argument as in [16, Lemma 3.4]). □

*Remark* 3.9. The preceeding lemma implies that the time derivative of the approximation of the pressure $\partial_t p_\varepsilon = \partial_t |n_h|^\gamma = g_\varepsilon + k_\varepsilon$ where $g_\varepsilon$ is uniformly bounded in $L^1([0,T] \times \Omega)$ and $k_\varepsilon$ in $L^1([0,T]; H^{-1}(\Omega))$. Hence $\partial_t W_\varepsilon = U_\varepsilon + V_\varepsilon$ where $U_\varepsilon \in L^1([0,T]; H^1(\Omega))$ solves $-\mu \Delta U_\varepsilon + U_\varepsilon = k_\varepsilon$ and $V_\varepsilon \in L^1([0,T]; W^{1,r}(\Omega))$, $1 \leq r < 1^* := d/(d-1)$ solves $-\mu \Delta V_\varepsilon + V_\varepsilon = g_\varepsilon$ (see [2, Thm. 6.1] for a proof of the second statement). Hence $\partial_t W_\varepsilon \in L^1([0,T]; W^{1,r}(\Omega))$ for any $1 \leq r < 1^*$.

### 3.3. Passing to the limit $\varepsilon \to 0$.

The estimates of the previous (sub)sections allow us to pass to the limit $\varepsilon \to 0$ in a subsequence still denoted $\varepsilon$ and conclude existence of limit functions

$$n_\varepsilon \rightharpoonup n \geq 0, \quad \text{in } L^q([0,T] \times \Omega), 1 \leq q < \infty,$$
$$p_\varepsilon \rightharpoonup \overline{p} \geq 0, \quad \text{in } L^q([0,T] \times \Omega), 1 \leq q < \infty,$$

where $p_\varepsilon := n_\varepsilon^\gamma$ and $0 \leq n, \overline{p} \in L^\infty([0,T] \times \Omega)$. Using Aubin-Lions' lemma for $W_\varepsilon, c_\varepsilon, q_\varepsilon$ and $\nabla W_\varepsilon$, we obtain strong convergence of a subsequence in $L^q([0,T] \times \Omega)$ for any $q \in [0, \infty)$ to limit functions $W, \nabla W \in L^q([0,T] \times \Omega)$ and $c, q \in L^2([0,T]; H^1(\Omega))$. Moreover, from the estimates in Lemma 3.7 we obtain that $W \in L^\infty([0,T] \times \Omega) \cap L^\infty([0,T]; H^2(\Omega))$. Hence we have that $(n, W, \overline{p}, c, q)$ satisfy for any $\varphi, \psi_i \in C_0^1([0,T) \times \Omega)$, $i = 1, 2, 3$,

$$\int_0^T \int_\Omega n \varphi_t - n \nabla W \cdot \nabla \varphi \, dx dt + \int_\Omega n_0 \, \varphi(0, x) dx = -\int_0^T \int_\Omega \overline{n \mathbf{\Phi}(p, c, q)} \varphi \, dx dt$$

$$\int_0^T \int_\Omega W \psi_1 + \mu \nabla W \cdot \nabla \psi_1 \, dx dt = \int_0^T \int_\Omega \overline{p} \, \psi_1 \, dx dt$$

$$\int_0^T \int_\Omega c \, \partial_t \psi_2 - \nu_c \nabla c \cdot \nabla \psi_2 \, dx dt = -\int_0^T \int_\Omega \left( c \overline{\mathbf{\Psi}_c(n, c)} + r_c(c_{\mathrm{supp}} - c) \right) \psi_2 \, dx dt$$

$$\int_0^T \int_\Omega q \, \partial_t \psi_3 - \nu_q \nabla q \cdot \nabla \psi_3 \, dx dt = -\int_0^T \int_\Omega \left( q \overline{\mathbf{\Psi}_q(n, q)} + r_q(q_{\mathrm{supp}} - q) \right) \psi_3 \, dx dt$$
(3.4)

where $\overline{n\boldsymbol{\Phi}(p,c,q)}$ is the weak limit of $n_\varepsilon\boldsymbol{\Phi}(p_\varepsilon,c_\varepsilon,q_\varepsilon)$. To conclude that the limit $(n,W,p,c,q)$ is a weak solution of (1.9), we need to show that $n_\varepsilon$ converges strongly and therefore in the limit $\overline{p} = p := n^\gamma$ and $\overline{n\boldsymbol{\Phi}(p,c,q)} = n\boldsymbol{\Phi}(p,c,q)$, $\overline{\boldsymbol{\Psi}_c(n,c)} = \boldsymbol{\Psi}_c(n,c)$ and $\overline{\boldsymbol{\Psi}_q(n,q)} = \boldsymbol{\Psi}_q(n,q)$. For this purpose, we combine a compensated compactness property (Lemma 3.11) with a monotonicity argument. We will also make use of the following lemma which was proved in a more general form in [6, 13] and which we proved in this particular form in [16, Lemma 3.6]:

**Lemma 3.10.** *Let $n, f \in L^\infty([0,T] \times \Omega)$ and $\boldsymbol{u} \in L^\infty([0,T]; H^1(\Omega))$ with $\mathrm{div}\,\boldsymbol{u} \in L^\infty([0,T] \times \Omega)$ satisfy*

$$n_t - \mathrm{div}(\boldsymbol{u}n) = f, \tag{3.5}$$

*in the sense of distributions. Then they satisfy for all continuously differentiable functions $b \in C^1(\mathbb{R})$*

$$b(n)_t - \mathrm{div}(\boldsymbol{u}b(n)) = b'(n)f + [b'(n)n - b(n)]\,\mathrm{div}\,\boldsymbol{u}, \tag{3.6}$$

*in the sense of distribution.*

*Proof.* The proof of this Lemma can be found in the appendix. $\square$

Applying Lemma 3.10 for the weak limit $n$ in (3.4) with $b(n) = n^2$, we obtain that $n$ satisfies

$$\int_0^T \int_\Omega n^2\varphi_t - n^2\nabla W \cdot \nabla\varphi\,dxdt = -\int_0^T \int_\Omega (2n\overline{n\boldsymbol{\Phi}(p,c,q)} + n^2\Delta W)\varphi\,dxdt \tag{3.7}$$

for any test functions $\varphi \in C_0^1((0,T) \times \Omega)$. Besides that, we obtain from integrating (3.2) for $f(n) = n^2$ in space and time

$$\int_\Omega n_\varepsilon^2(\tau)\,dx - \int_\Omega n_\varepsilon^2(0)\,dx \leq \int_0^\tau \int_\Omega n_\varepsilon^2 \Delta W_\varepsilon + 2n_\varepsilon^2 \boldsymbol{\Phi}(p_\varepsilon, c_\varepsilon, q_\varepsilon)\,dxdt$$

Passing to the limit $\varepsilon \to 0$ in this inequality, we have

$$\int_\Omega \overline{n^2}(\tau)\,dx - \int_\Omega n_0^2\,dx \leq \int_0^\tau \int_\Omega \overline{n^2 \Delta W} + 2\overline{n_\varepsilon^2 \boldsymbol{\Phi}(p,c,q)}\,dxdt, \tag{3.8}$$

where $\overline{n^2}$ denotes the weak limit of $n_\varepsilon^2$ and $\overline{n^2 \Delta W}$ and $\overline{n^2\boldsymbol{\Phi}(p,c,q)}$ are the weak limits of $n_\varepsilon^2 \Delta W_\varepsilon$ and $n_\varepsilon^2 \boldsymbol{\Phi}(p_\varepsilon, c_\varepsilon, q_\varepsilon)$ respectively. Letting $\tau \to 0$ in this inequality, we obtain, thanks to the boundedness of the integrand on the right hand side,

$$\int_\Omega \overline{n^2}(0)\,dx - \int_\Omega n_0^2\,dx \leq 0.$$

Therefore, since $b(n) = n^2$ is convex, $\overline{n^2} \geq n^2$ and so $\overline{n^2}(0,x) = n_0^2(x)$. We now choose smooth test functions $\varphi_\epsilon$ approximating $\varphi(t,x) = \mathbf{1}_{[0,\tau]}(t)$, where $\tau \in (0,T]$, in inequality (3.7) and then pass to the limit in the approximation to obtain the inequality

$$\int_\Omega n^2(\tau)\,dx - \int_\Omega n_0^2\,dx = \int_0^\tau \int_\Omega (2n\overline{n\boldsymbol{\Phi}(p,c,q)} + n^2\Delta W)\,dxdt \tag{3.9}$$

Subtracting (3.9) from (3.8), we have

$$\int_\Omega \left(\overline{n^2} - n^2\right)(\tau)dx$$
$$\leq \int_0^\tau \int_\Omega \left(2\overline{n^2\boldsymbol{\Phi}(p,c,q)} - 2n\overline{n\boldsymbol{\Phi}(p,c,q)} + \Delta W\left(\overline{n^2} - n^2\right) + \overline{n^2\Delta W} - \overline{n^2}\Delta W\right)dx\,dt. \tag{3.10}$$

Now using the explicit expression of $\boldsymbol{\Phi}$, (1.2), the first term on the right hand side can be estimated as follows:

$$\int_0^\tau \int_\Omega \left(2\overline{n^2 \boldsymbol{\Phi}(p,c,q)} - 2n\overline{n\boldsymbol{\Phi}(p,c,q)}\right) dx\, dt$$
$$= 2\int_0^\tau \int_\Omega \left(\overline{n^2 g_1(c,q)\boldsymbol{G}(p)} - \overline{n^2 g_2(c,q)} - n\overline{ng_1(c,q)\boldsymbol{G}(p)} + n\overline{ng_2(c,q)}\right) dx\, dt$$
$$= 2\int_0^\tau \int_\Omega \left(g_1(c,q)\left(\overline{n^2 \boldsymbol{G}(p)} - n\overline{n\boldsymbol{G}(p)}\right) - g_2(c,q)\left(\overline{n^2} - n^2\right)\right) dx\, dt$$
$$\leq 2\int_0^\tau \int_\Omega g_1(c,q)\left(\overline{n^2 \boldsymbol{G}(p)} - n\overline{n\boldsymbol{G}(p)}\right) dx\, dt$$
$$= 2\int_0^\tau \int_\Omega \alpha\, g_1(c,q)\left(\overline{n^2} - n^2\right) - \beta\, g_1(c,q)\left(\overline{n^2 g_0(n)} - n\overline{ng_0(n)}\right) dx\, dt$$
$$\leq 2\int_0^\tau \int_\Omega \alpha\, g_1(c,q)\left(\overline{n^2} - n^2\right) - \beta\, g_1(c,q)\left(\overline{n^2 g_0(n)} - \overline{n^2} g_0(n)\right) dx\, dt$$
$$\leq 2\alpha \|g_1\|_\infty \int_0^\tau \int_\Omega \left(\overline{n^2} - n^2\right) dx\, dt \tag{3.11}$$

The second inequality follows from the strong convergence of $\{c_\varepsilon\}_{\varepsilon>0}$ and $\{q_\varepsilon\}_{\varepsilon>0}$, the first inequality from the nonnegativity of $g_2$ and the convexity of $b(n) = n^2$, for the second last inequality, we used [13, Lemma 3.35] and for the last inequality the boundedness of $g_1$. To estimate the second term on the right hand side of (3.10), we use that $\Delta W$ is bounded thanks to Lemma 3.7 and that $\overline{n^2} \geq n^2$ by the convexity of $f(x) = x^2$. Hence

$$\int_0^\tau \int_\Omega \Delta W \left(\overline{n^2} - n^2\right) dx dt \leq \frac{P_M}{\mu} \int_0^\tau \int_\Omega \left(\overline{n^2} - n^2\right) dx dt. \tag{3.12}$$

In order to deal with the last term, we use the following lemma from (cf. Lemma 3.7, Trivisa and Weber [16]). We present here the main steps of the proof for completeness.

**Lemma 3.11.** *The weak limits $(n, W, \overline{p})$ of the sequences $\{(n_\varepsilon, W_\varepsilon, p_\varepsilon)\}_{\varepsilon>0}$ satisfy for smooth functions $S : \mathbb{R} \to \mathbb{R}$,*

$$\int_\Omega \left(\overline{S(n)\Delta W} - \overline{S(n)}\Delta W\right) dx = \frac{1}{\mu} \int_\Omega \left(\overline{p}\,\overline{S(n)} - \overline{pS(n)}\right) dx \tag{3.13}$$

*where $\overline{S(n)\Delta W}$, $\overline{S(n)}$, $\overline{pS(n)}$ are the weak limits of $S(n_\varepsilon)\Delta W_\varepsilon$, $S(n_\varepsilon)$ and $p_\varepsilon S(n_\varepsilon)$ respectively.*

*Proof.* We multiply the equation for $W_\varepsilon$ by $S(n_\varepsilon)$ and integrate over $\Omega$,

$$\int_\Omega \mu \Delta W_\varepsilon\, S(n_\varepsilon) - W_\varepsilon S(n_\varepsilon)\, dx = -\int_\Omega p_\varepsilon S(n_\varepsilon)\, dx.$$

Passing to the limit $\varepsilon \to 0$, we obtain

$$\int_\Omega \mu \overline{\Delta W S(n)} - W \overline{S(n)}\, dx = -\int_\Omega \overline{pS(n)}\, dx. \tag{3.14}$$

On the other hand, using the smooth function $S(n_\varepsilon)$ as a test function in the weak formulation of the limit equation

$$-\mu \Delta W + W = \overline{p},$$

and passing to the limit $\varepsilon \to 0$, we obtain

$$\int_\Omega \mu \Delta W \overline{S(n)} - W \overline{S(n)}\, dx = -\int_\Omega \overline{p}\,\overline{S(n)}\, dx.$$

Combining the last identity with (3.14), we obtain (3.13). □

Applying this lemma to the second term in (3.10) with $S(n) = n^2$, we can estimate it by

$$\int_0^\tau \int_\Omega \left( \overline{n^2 \Delta W} - \overline{n^2} \Delta W \right) dx = \frac{1}{\mu} \int_\Omega \left( \overline{p}\,\overline{n^2} - \overline{pn^2} \right) dx dt$$
$$= \frac{1}{\mu} \int_\Omega \left( \overline{n^\gamma}\,\overline{n^2} - \overline{n^{2+\gamma}} \right) dx dt \leq 0,$$

using that $\overline{n^\gamma}\,\overline{n^2} \leq \overline{n^{2+\gamma}}$ (cf. [13]). Thus,

$$\int_\Omega \left( \overline{n^2} - n^2 \right)(\tau) dx \leq \left( 2\alpha \|g_1\|_\infty + \frac{P_M}{\mu} \right) \int_0^\tau \int_\Omega \left( \overline{n^2} - n^2 \right) dx\, dt.$$

Hence Grönwall's inequality implies $\int_\Omega \left( \overline{n^2} - n^2 \right)(\tau) dx \leq 0$. By convexity of the function $b(n) = n^2$ we have $n^2 \leq \overline{n^2}$ almost everywhere and hence $\overline{n^2}(t, x) = n^2(t, x)$ almost everywhere in $(0, T) \times \Omega$. Therefore we conclude that the functions $n_\varepsilon$ converge strongly to $n$ almost everywhere and in particular also $\bar{p} = n^\gamma$ which means that the limit $(n, W, \bar{p}, c, q)$ is a weak solution of the equations (1.9).

## 4. Global existence via a numerical approximation

In this section, we will construct a finite difference scheme to approximate (1.9) on a rectangular domain, for simplicity we will use $\Omega = [0, 1]^d$ where $d = 1, 2, 3$ is the spatial dimension; the extension of the scheme to other rectangular domains with nonuniform mesh widths is straighforward. The assumptions on initial and boundary data are the same as in Section 3. We denote $N_x \in \mathbb{N}$ the number of grid cells in one coordinate direction, $h := 1/N_x$ the mesh width and $N_t = N_x/\kappa \in \mathbb{N}$ the number of time steps and $\Delta t := 1/N_t$ the time step size. We will determine conditions on the time step size $\Delta t = \kappa h > 0$ later on. We define gridpoints and grid cells

$$\mathcal{C}_{i_1,\ldots,i_d} := ((i_1 - 1)h, i_1 h] \times \cdots \times ((i_d - 1)h, i_d h],$$
$$x_{i_1,\ldots,i_d} := ((i_1 - 1/2)h, \ldots, (i_d - 1/2)h),$$

and time steps $t^m := m\Delta t$, $m = 0, \ldots, N_t$. To simplify notation, we introduce the multiindex $\underline{i} \in \mathcal{I}_{N_x} := \{0, \ldots, N_x\}^d$, such that $\underline{i} = (i_1, \ldots, i_d)$ and we can write

$$\mathcal{C}_{\underline{i}} = \mathcal{C}_{i_1,\ldots,i_d}, \quad x_{\underline{i}} = x_{i_1,\ldots,i_d}.$$

We will approximate $(n, W, p, c, q)$ at these points. Specifically,

$$f_{\underline{i}}^m \approx f(m\Delta t, x_{\underline{i}}).$$

where $f \in \{n, W, p, c, q\}$. Next, we let $\mathbf{e}_1 := (1, 0, 0)$, $\mathbf{e}_2 := (0, 1, 0)$, and $\mathbf{e}_3 := (0, 0, 1)$. Using these vectors, we then define the forward and backward difference operators

$$D_j^+ f_{\underline{i}} = \frac{f_{\underline{i}+\mathbf{e}_j} - f_{\underline{i}}}{h}, \qquad D_j^- f_{\underline{i}} = D_j^+ f_{\underline{i}-\mathbf{e}_j}, \qquad D_t^\pm f^m = \pm \frac{f^{m\pm 1} - f^m}{\Delta t}.$$

respectively, for $j = 1, \ldots, d$, and $\underline{i} \in \mathcal{I}_{N_x}$. Based on these, we define the discrete Laplace, divergence and gradient operators,

$$\nabla_h^\pm := (D_1^\pm, \ldots, D_d^\pm)^T, \quad \operatorname{div}_h^\pm v_{\underline{i}} := \sum_{j=1}^d D_j^\pm v^{(j)}, \quad \Delta_h := \operatorname{div}_h^\pm \nabla_h^\mp.$$

Whenever the choice of the forward $\nabla_h^+$ or backward $\nabla_h^-$ difference operator does not matter, we will write $\nabla_h$ or $\operatorname{div}_h$ respectively.

4.1. **An explicit finite difference scheme.** Given $(n_{\underline{i}}^m, W_{\underline{i}}^m, c_{\underline{i}}^m, q_{\underline{i}}^m)$ at time step $m$, we define the quantities $(n_{\underline{i}}^{m+1}, W_{\underline{i}}^{m+1}, c_{\underline{i}}^{m+1}, q_{\underline{i}}^{m+1})$ at the next time step by

$$-\mu \Delta_h W_{\underline{i}}^m + W_{\underline{i}}^m = p_{\underline{i}}^m := |n_{\underline{i}}^m|^\gamma, \tag{4.1a}$$

$$D_t^+ n_{\underline{i}}^m + \operatorname{div}_h^- \mathbf{F}_{\underline{i}}(W_{\underline{i}}^m, n_{\underline{i}}^m) = n_{\underline{i}}^m \boldsymbol{\Phi}(p_{\underline{i}}^m, c_{\underline{i}}^m, q_{\underline{i}}^m), \tag{4.1b}$$

where the flux $\mathbf{F}_{\underline{i}} = (F_{\underline{i}}^{(1)}, \ldots, F_{\underline{i}}^{(d)})$ is defined by

$$F_{\underline{i}}^{(j)}(W_{\underline{i}}, n_{\underline{i}}) = -D_j^+ W_{\underline{i}} \frac{n_{\underline{i}} + n_{\underline{i}+\mathbf{e}_j}}{2} - \frac{h}{2} |D_j^+ W_{\underline{i}}| D_j^+ n_{\underline{i}}. \tag{4.2}$$

The diffusion equations (1.9c) and (1.9d) naturally require a quadratic CFL-condition whereas the transport equation only requires a linear one. We therefore use operator splitting to solve the whole system, specifically, we evolve $n_{\underline{i}}^m$ with the larger time step $\Delta t$ and then compute $c_{\underline{i}}^{m+1}$ and $q_{\underline{i}}^{m+1}$ iterating the following scheme with time steps $\Delta t_k$, $k \in \{c, q\}$, satisfying (4.11), for $s = 1, \ldots, N_k := \operatorname{rnd}(\Delta t / \Delta t_k)$ times:

$$D_{t_c}^+ c_{\underline{i}}^{m,s} - \nu_c \Delta_h c_{\underline{i}}^{m,s} = c_{\underline{i}}^{m,s} \boldsymbol{\Psi}_c(n_{\underline{i}}^m, c_{\underline{i}}^{m,s}) + r_c(c_{\text{supp}} - c_{\underline{i}}^{m,s}), \tag{4.3a}$$

$$D_{t_q}^+ q_{\underline{i}}^{m,s} - \nu_q \Delta_h q_{\underline{i}}^{m,s} = q_{\underline{i}}^{m,s} \boldsymbol{\Psi}_q(n_{\underline{i}}^m, q_{\underline{i}}^{m,s}) + r_q(q_{\text{supp}} - q_{\underline{i}}^{m,s}), \tag{4.3b}$$

where $D_{t_k}^\pm f^{m,s} = \pm \frac{f^{m,s+1} - f^{m,s}}{\Delta t_k}$. We set $c_{\underline{i}}^{m,0} := c_{\underline{i}}^m$ and $q_{\underline{i}}^{m,0} := q_{\underline{i}}^m$, and $c_{\underline{i}}^{m,N_{t_c}} =: c_{\underline{i}}^{m+1}$ and $q_{\underline{i}}^{m,N_{t_q}} =: q_{\underline{i}}^{m+1}$.

We use homogeneous Neumann or periodic boundary conditions for the variables $n$ and $W$ ($\underline{i} = (i_1, \ldots, i_d)$):

$$\begin{aligned} n_{\underline{i}-\mathbf{e}_j}^m &= n_{\underline{i}}^m, \quad i_j = 1, & n_{\underline{i}+\mathbf{e}_j}^m &= n_{\underline{i}}^m, \quad i_j = N_x, \\ W_{\underline{i}-\mathbf{e}_j}^m &= W_{\underline{i}}^m, \quad i_j = 1, & W_{\underline{i}+\mathbf{e}_j}^m &= W_{\underline{i}}^m, \quad i_j = N_x, \end{aligned}$$

For the variables $c$ and $q$ we use the Dirichlet boundary conditions,

$$\begin{aligned} c_{\underline{i}-\mathbf{e}_j}^{m,s} &= c_b(m\Delta t + s\Delta t_c, x_{\underline{i}-\mathbf{e}_j}), \quad i_j = 1; & c_{\underline{i}+\mathbf{e}_j}^{m,s} &= c_b(m\Delta t + s\Delta t_c, x_{\underline{i}+\mathbf{e}_j}), \quad i_j = N_x; \\ q_{\underline{i}-\mathbf{e}_j}^{m,s} &= q_b(m\Delta t + s\Delta t_c, x_{\underline{i}-\mathbf{e}_j}), \quad i_j = 1; & q_{\underline{i}+\mathbf{e}_j}^{m,s} &= q_b(m\Delta t + s\Delta t_c, x_{\underline{i}+\mathbf{e}_j}), \quad i_j = N_x; \end{aligned}$$

$j = 1, \ldots, d$, $m = 1, \ldots, N_t$, $s = 1, \ldots, N_k$ $k \in \{c, q\}$. The initial condition we approximate taking averages over the cells,

$$n_{\underline{i}}^0 = \frac{1}{|\mathcal{C}_{\underline{i}}|} \int_{\mathcal{C}_{\underline{i}}} n_0(x) \, dx, \quad p_{\underline{i}}^0 = |n_{\underline{i}}^0|^\gamma, \quad \underline{i} \in \mathcal{I}_{N_x},$$

and analoguously we approximate $c_0$ and $q_0$.

*Remark* 4.1 (Implicit time stepping). Another option would be to compute the solutions to the diffusion equations (1.9c) and (1.9d) implicitly, that would improve the CFL-condition, however becomes expensive in several space dimension because matrices corresponding to the difference operators have to be inverted in every timestep.

4.2. **Estimates on approximations.** In the following, we will prove estimates on the discrete quantities $(n_{\underline{i}}^m, W_{\underline{i}}^m, c_{\underline{i}}^m, q_{\underline{i}}^m)$ obtained using the scheme (4.1)–(4.3). We therefore define the piecewise constant functions

$$f_h(t, x) = \sum_{m=0}^{N_t - 1} \sum_{\underline{i} \in \mathcal{I}_{N_x}} f_{\underline{i}}^m \mathbf{1}_{\mathcal{C}_{\underline{i}}}(x) \mathbf{1}_{[t^m, t^{m+1}]}(t), \quad (t, x) \in [0, T] \times \Omega, \tag{4.4}$$

where $f \in \{n, W, p\}$. For the nutrient $c$ and the drug $q$ we define

$$c_h(t,x) = \sum_{m=0}^{N_t-1} \sum_{s=0}^{N_c-1} \sum_{\underline{i} \in \mathcal{I}_{N_x}} c_{\underline{i}}^{m,s} \mathbf{1}_{\mathcal{C}_{\underline{i}}}(x) \mathbf{1}_{[t^m + s\Delta t_c, t^m + (s+1)\Delta t_c)}(t),$$

$$q_h(t,x) = \sum_{m=0}^{N_t-1} \sum_{s=0}^{N_q-1} \sum_{\underline{i} \in \mathcal{I}_{N_x}} q_{\underline{i}}^{m,s} \mathbf{1}_{\mathcal{C}_{\underline{i}}}(x) \mathbf{1}_{[t^m + s\Delta t_q, t^m + (s+1)\Delta t_q)}(t),$$

for $(t,x) \in [0,T] \times \Omega$. We first prove that $n_h$ stays nonnegative and uniformly bounded from above.

**Lemma 4.2.** *For any $t > 0$, the functions $n_h(t,\cdot)$ are uniformly (in $h > 0$) bounded and nonnegative if $0 \leq n_{\underline{i}}^0 \leq n_\infty := P_M^{1/\gamma} < \infty$ uniformly in $h > 0$ and the timestep $\Delta t$ satisfies the CFL condition*

$$\Delta t \leq \min\left\{\frac{h}{4d\max_{\underline{i},j}|D_j^+ W_{\underline{i}}^m| + h\mathbf{\Phi}^\infty}, \frac{\mu}{4\gamma \overline{n}_\infty^\gamma}\right\} \tag{4.5}$$

*(where $\mathbf{\Phi}^\infty := \sup_{0 \leq s_i \leq s_{i,\infty}} \mathbf{\Phi}(s_1, s_2, s_3)$, $(s_{1,\infty} := \infty, s_{2,\infty} := c_\infty, s_{3,\infty} := q_\infty))$*

$$\overline{n}_\infty = n_\infty + 4\Delta t \sup_{0 \leq s_i \leq s_{i,\infty}} \left(s_1^{1/\gamma} \mathbf{\Phi}(s_1, s_2, s_3)\right),$$

*we have for all $m \geq 0$,*

$$0 \leq \min_{\underline{i}} n_{\underline{i}}^m \leq \max_{\underline{i}} n_{\underline{i}}^m \leq \overline{n}_\infty. \tag{4.6}$$

*Proof.* The proof goes by induction on the timestep $m$. Clearly, by the assumptions, we have $0 \leq n_{\underline{i}}^0 \leq \overline{n}_\infty$. For the induction step we therefore assume that this holds for timestep $m \geq 0$ and show that it implies the nonnegativity and boundedness at timestep $m+1$.

We first show that the $W_{\underline{i}}^m$ are bounded in terms of the $p_{\underline{i}}^m$. To do so, let us assume it has a local maximum $W_{\underline{\hat{i}}}^m$ in a cell $\mathcal{C}_{\underline{\hat{i}}}$, for some $\underline{\hat{i}} \in \mathcal{I}_{N_x}$. Then

$$D_j^+ W_{\underline{\hat{i}}}^m \leq 0, \quad -D_j^- W_{\underline{\hat{i}}}^m \leq 0, \quad j = 1, \ldots, d,$$

(if $\hat{i}_j \in \{1, N_x\}$ for some $j \in \{1, \ldots, d\}$, then because of the Neumann boundary conditions, the forward/backward difference in direction of the boundary is zero and thus the previous inequality is true as well). Hence

$$\Delta_h W_{\underline{\hat{i}}}^m = \frac{1}{h} \sum_{j=1}^d \left(D_j^+ W_{\underline{\hat{i}}}^m - D_j^- W_{\underline{\hat{i}}}^m\right) \leq 0.$$

Therefore,

$$W_{\underline{\hat{i}}}^m = p_{\underline{\hat{i}}}^m + \mu \Delta_h W_{\underline{\hat{i}}}^m \leq p_{\underline{\hat{i}}}^m \leq \max_{\underline{i}} |n_{\underline{i}}^m|^\gamma.$$

Similarly, at a local minimum $W_{\underline{\hat{i}}}^m$ of $W_h$, we have

$$D_j^+ W_{\underline{\hat{i}}}^m \geq 0, \quad -D_j^- W_{\underline{\hat{i}}}^m \geq 0, \quad j = 1, \ldots, d,$$

and hence

$$\Delta_h W_{\underline{\hat{i}}}^m = \frac{1}{h} \sum_{j=1}^d \left(D_j^+ W_{\underline{\hat{i}}}^m - D_j^- W_{\underline{\hat{i}}}^m\right) \geq 0,$$

which implies

$$W_{\underline{\hat{i}}}^m = p_{\underline{\hat{i}}}^m + \mu \Delta_h W_{\underline{\hat{i}}}^m \geq p_{\underline{\hat{i}}}^m \geq \min_{\underline{i}} |n_{\underline{i}}^m|^\gamma \geq 0.$$

Thus,

$$0 \leq W_h \leq \max_{\underline{i}} |n_{\underline{i}}^m|^\gamma. \tag{4.7}$$

Now we rewrite the scheme (4.1b) as

$$n_{\underline{i}}^{m+1} = \left(\alpha_{\underline{i}}^{(1),m} + \alpha_{\underline{i}}^{(2),m}\right) n_{\underline{i}}^m + \sum_{j=1}^d \left(\beta_{\underline{i},j}^{m,+} n_{\underline{i}+\mathbf{e}_j}^m + \beta_{\underline{i},j}^{m,-} n_{\underline{i}-\mathbf{e}_j}^m\right) \quad (4.8)$$

where

$$\alpha_{\underline{i}}^{(1),m} = 1 - \frac{\Delta t}{2}\Delta_h W_{\underline{i}}^m - \frac{\Delta t}{2h}\sum_{j=1}^d \left(|D_j^+ W_{\underline{i}}^m| + |D_j^- W_{\underline{i}}^m|\right)$$

$$\alpha_{\underline{i}}^{(2),m} = \Delta t\, \boldsymbol{\Phi}(p_{\underline{i}}^m, c_{\underline{i}}^m, q_{\underline{i}}^m) + \Delta t \Delta_h W_{\underline{i}}^m,$$

$$\beta_{\underline{i},j}^{m,\pm} = \frac{\Delta t}{2h}\left(|D_j^\pm W_{\underline{i}}^m| \pm D_j^\pm W_{\underline{i}}^m\right), \quad j = 1, \ldots, d.$$

We note that $\beta_{\underline{i},j}^{m,\pm} \geq 0$, and that under the CFL-condition (4.5), also $\alpha_{\underline{i}}^{(1),m} + \alpha_{\underline{i}}^{(2),m} \geq 0$. Therefore, under the assumption that $n_{\underline{i}}^m \geq 0$ for all $\underline{i} \in \mathcal{I}_{N_x}$, we have

$$n_{\underline{i}}^{m+1} \geq \sum_{j=1}^d \left(\beta_{\underline{i},j}^{m,+} + \beta_{\underline{i},j}^{m,-}\right) \min\left\{\min_j n_{\underline{i}+\mathbf{e}_j}^m, \min_j n_{\underline{i}-\mathbf{e}_j}^m\right\}$$
$$+ \left(\alpha_{\underline{i}}^{(1),m} + \alpha_{\underline{i}}^{(2),m}\right) n_{\underline{i}}^m \geq 0.$$

We proceed to showing the boundedness of $n_h$. Thanks to the CFL-condition (4.5), we have

$$\alpha_{\underline{i}}^{(1),m} \geq \frac{1}{2}, \quad \beta_{\underline{i},j}^{m,\pm} \leq \frac{1}{4d}.$$

Furthermore, $\alpha_{\underline{i}}^{(1),m} + \sum_{j=1}^d (\beta_{\underline{i},j}^{m,+} + \beta_{\underline{i},j}^{m,-}) = 1$. Using the induction hypothesis that $n_{\underline{i}}^m \leq \overline{n}_\infty$ for all $\underline{i}$ and the nonnegativity of $n_h$ which we have just proved, we can estimate $n_{\underline{i}}^{m+1}$:

$$n_{\underline{i}}^{m+1} \leq \left(\alpha_{\underline{i}}^{(1),m} + \alpha_{\underline{i}}^{(2),m}\right) n_{\underline{i}}^m + \sum_{j=1}^d \left(\beta_{\underline{i},j}^{m,+} + \beta_{\underline{i},j}^{m,-}\right) \overline{n}_\infty$$
$$\leq \left(\frac{1}{2} + \alpha_{\underline{i}}^{(2),m}\right) n_{\underline{i}}^m + \frac{1}{2}\overline{n}_\infty = \overline{n}_\infty - \frac{1}{2}\left(\overline{n}_\infty - n_{\underline{i}}^m\right) + \alpha_{\underline{i}}^{(2),m} n_{\underline{i}}^m \quad (4.9)$$

We can rewrite and bound $\alpha_{\underline{i}}^{(2),m}$ using the equation for $W_{\underline{i}}^m$, (4.1a),

$$\alpha_{\underline{i}}^{(2),m} = \Delta t \left(\boldsymbol{\Phi}(p_{\underline{i}}^m, c_{\underline{i}}^m, q_{\underline{i}}^m) + \Delta_h W_{\underline{i}}^m\right)$$
$$= \Delta t \left(\boldsymbol{\Phi}(p_{\underline{i}}^m, c_{\underline{i}}^m, q_{\underline{i}}^m) + \frac{1}{\mu}\left(W_{\underline{i}}^m - p_{\underline{i}}^m\right)\right)$$
$$\leq \Delta t \left(\boldsymbol{\Phi}(p_{\underline{i}}^m, c_{\underline{i}}^m, q_{\underline{i}}^m) + \frac{1}{\mu}\left(\overline{n}_\infty^\gamma - |n_{\underline{i}}^m|^\gamma\right)\right)$$
$$\leq \Delta t \left(\boldsymbol{\Phi}(p_{\underline{i}}^m, c_{\underline{i}}^m, q_{\underline{i}}^m) + \frac{\gamma\,\overline{n}_\infty^{\gamma-1}}{\mu}\left(\overline{n}_\infty - n_{\underline{i}}^m\right)\right)$$
$$\leq \Delta t\, \boldsymbol{\Phi}(p_{\underline{i}}^m, c_{\underline{i}}^m, q_{\underline{i}}^m) + \frac{1}{4\overline{n}_\infty}(\overline{n}_\infty - n_{\underline{i}}^m),$$

where we have used (4.7) for the first inequality, that $f(a) - f(b) = f'(\widetilde{a})(a-b)$ for some intermediate value $\widetilde{a} \in [b, a]$, with $f(a) = a^\gamma$, for the second inequality and the CFL-condition for the last inequality. Now going back to (4.9) and inserting this there, we obtain,

$$n_{\underline{i}}^{m+1} \leq \overline{n}_\infty - \frac{1}{2}\left(\overline{n}_\infty - n_{\underline{i}}^m\right) + \left(\Delta t \boldsymbol{\Phi}(p_{\underline{i}}^m, c_{\underline{i}}^m, q_{\underline{i}}^m) + \frac{1}{4\overline{n}_\infty}(\overline{n}_\infty - n_{\underline{i}}^m)\right) n_{\underline{i}}^m$$

$$\leq \frac{3}{4}\overline{n}_\infty + \frac{1}{4}n_{\underline{i}}^m + \Delta t\, n_{\underline{i}}^m\, \mathbf{\Phi}(p_{\underline{i}}^m, c_{\underline{i}}^m, q_{\underline{i}}^m) \tag{4.10}$$

If $n_{\underline{i}}^m \geq n_\infty$ then $\mathbf{\Phi}(p_{\underline{i}}^m, c_{\underline{i}}^m, q_{\underline{i}}^m) \leq 0$ and hence the expression in (4.10) is bounded by $\overline{n}_\infty$. On the other hand, if $n_{\underline{i}}^m \leq n_\infty$, we can bound it by

$$\begin{aligned} n_{\underline{i}}^{m+1} &\leq \frac{3}{4}\overline{n}_\infty + \frac{1}{4}n_{\underline{i}}^m + \Delta t\, n_{\underline{i}}^m\, \mathbf{\Phi}(p_{\underline{i}}^m, c_{\underline{i}}^m, q_{\underline{i}}^m) \\ &\leq \frac{3}{4}\overline{n}_\infty + \frac{1}{4}\left(n_\infty + 4\Delta t \sup_{0 \leq s_i \leq s_{i,\infty}} \left(s_1^{1/\gamma}\mathbf{\Phi}(s_1, s_2, s_3)\right)\right) = \overline{n}_\infty \end{aligned}$$

where we used the definition of $\overline{n}_\infty$ for the last equality. This proves that $n_{\underline{i}}^{m+1} \leq \overline{n}_\infty$ for all $\underline{i}$ if the same holds already for the $n_{\underline{i}}^m$. $\square$

*Remark* 4.3. The estimates in the proof of the previous lemma are quite coarse and balancing the terms in a more carefu way one ends up with a better CFL-condition than (4.5) in practice. Also note that $\overline{n}_\infty \to n_\infty$ when $\Delta t \to 0$.

Next we prove that the approximations $\{c_h\}_{h>0}$ of the nutrient function are nonnegative and uniformly bounded. A similar argument can be used to show that the sequence of approximations $\{q_h\}_{h>0}$ is nonnegative and uniformly bounded.

**Lemma 4.4.** *If*

$$\Delta t_c \leq \min\left\{\frac{h\Delta t}{\nu_c}, \frac{h^2}{2^d\nu_c + h^2 r_c + h^2 \mathbf{\Psi}_c^\infty}\right\}, \tag{4.11}$$

*where $\Delta t > 0$ satisfies (4.5), and*

$$\mathbf{\Psi}_c^\infty := \max_{\substack{0 \leq n \leq \overline{n}_\infty \\ 0 \leq c \leq c_\infty}} |\mathbf{\Psi}_c(n, c)|,$$

*and the initial data $c_0$ and its approximation $c_{\underline{i}}^0$, as well as $c_b$ and $c_{supp}$, are nonnegative and uniformly bounded, that is,*

$$0 \leq c_{\underline{i}}^0, c_b, c_{supp} \leq c_\infty < \infty,$$

*then we have uniformly in $h > 0$, for all $m \in \mathbb{N}$, $s = 0, \ldots, N_c$,*

$$0 \leq \min_{\underline{i}} c_{\underline{i}}^{m,s} \leq \max_{\underline{i}} c_{\underline{i}}^{m,s} \leq c_\infty.$$

*Proof.* We prove the lemma by induction on $s, m$. We assume that the claim holds for the approximations $c_{\underline{i}}^{m,s}$ and show that this implies that it holds for $c_{\underline{i}}^{m,s+1}$ as well. To simplify the notation, we omit writing the index $m$. We start by proving the nonnegativity. Again, we rewrite the scheme (as we did in Lemma 4.2) as

$$c_{\underline{i}}^{s+1} = \left(\alpha_{\underline{i}}^{(1),s} + \alpha_{\underline{i}}^{(2),s}\right) c_{\underline{i}}^s + \sum_{j=1}^d \left(\beta_{\underline{i},j}^{s,+} c_{\underline{i}+\mathbf{e}_j}^s + \beta_{\underline{i},j}^{s,-} c_{\underline{i}-\mathbf{e}_j}^s\right) + \gamma \tag{4.12}$$

where

$$\begin{aligned} \alpha_{\underline{i}}^{(1),s} &= 1 - \frac{2^d \nu_c \Delta t_c}{h^2}, \quad \alpha_{\underline{i}}^{(2),s} = \Delta t_c \left(\mathbf{\Psi}_c(n_{\underline{i}}, c_{\underline{i}}^s) - r_c\right), \\ \beta_{\underline{i},j}^{s,\pm} &= \frac{\nu_c \Delta t_c}{h^2}, \quad j = 1, \ldots, d, \quad \gamma = \Delta t_c\, r_c\, c_{\text{supp}}. \end{aligned} \tag{4.13}$$

We note that under the CFL-condition (4.11), $(\alpha_{\underline{i}}^{(1),s} + \alpha_{\underline{i}}^{(2),s}), \beta_{\underline{i},j}^{s,\pm} \geq 0$. Furthermore, $\gamma \geq 0$. Therefore, since by assumption $c_{\underline{i}}^s \geq 0$ for all $\underline{i} \in \mathcal{I}_{N_x}$,

$$\begin{aligned} c_{\underline{i}}^{s+1} &\geq \sum_{j=1}^d \left(\beta_{\underline{i},j}^{s,+} + \beta_{\underline{i},j}^{s,-}\right) \min\left\{\min_j c_{\underline{i}+\mathbf{e}_j}^s, \min_j c_{\underline{i}-\mathbf{e}_j}^s\right\} \\ &\quad + \left(\alpha_{\underline{i}}^{(1),s} + \alpha_{\underline{i}}^{(2),s}\right) c_{\underline{i}}^s + \gamma \geq 0. \end{aligned} \tag{4.14}$$

To prove the uniform boundedness, we proceed in a very similar way. We define $\widehat{c}_{\underline{i}}^s := c_\infty - c_{\underline{i}}^s$ and show that $\widehat{c}_{\underline{i}}^{s+1}$ is nonnegative if $\widehat{c}_{\underline{i}}^s$ is nonnegative (that $\widehat{c}_{\underline{i}}^0$ follows from the assumptions). The scheme for $\widehat{c}_{\underline{i}}^m$ is the following:

$$D_{t_c}^+ \widehat{c}_{\underline{i}}^s - \nu_c \Delta_h \widehat{c}_{\underline{i}}^s = (\widehat{c}_{\underline{i}}^s - c_\infty)\boldsymbol{\Psi}_c(n_{\underline{i}}, c_\infty - \widehat{c}_{\underline{i}}^s) + r_c(c_\infty - c_{\text{supp}} - \widehat{c}_{\underline{i}}^s).$$

We observe that this is in fact (4.12) – (4.13) with $\alpha_{\underline{i}}^{(2),s}$ and $\gamma$ replaced by

$$\alpha_{\underline{i}}^{(2),s} = \Delta t_c \left(\boldsymbol{\Psi}_c(n_{\underline{i}}, c_\infty - c_{\underline{i}}^s) - r_c\right),$$
$$\gamma = \Delta t_c \left(r_c(c_\infty - c_{\text{supp}}) - c_\infty \boldsymbol{\Psi}_c(n_{\underline{i}}, c_\infty - c_{\underline{i}}^s)\right).$$

Thanks to the properties of $\boldsymbol{\Psi}_c$, $\gamma \geq 0$ and $\alpha_{\underline{i}}^{(1),s} + \alpha_{\underline{i}}^{(2),s} \geq 0$ and we conclude (4.14) for $\widehat{c}_{\underline{i}}^{s+1}$ which implies the desired maximum principle for $c_h$. $\square$

In analogy to the continuous case, we would like to prove estimates on higher order differences of $c_h$ to conclude strong convergence of the sequence $\{c_h\}_{h>0}$.

**Lemma 4.5.** *Let $\Delta t_c$ satisfy the CFL-condition (4.11), $\Delta t$ satisfy (4.5) and $0 \leq c_{\underline{i}}^0, c_b, c_{supp} \leq c_\infty$, $0 \leq n_{\underline{i}}^0 \leq n_\infty$, $c_b$ such that it has a unique extension $\widetilde{c}_b \in C^2([0,T] \times \overline{\Omega})$. Then $\nabla_h c_h \in L^2((0,T) \times \Omega)$ uniformly in $h > 0$ (where $\nabla_h c_h = \nabla_h^+ c_h$ or $\nabla_h c_h = \nabla_h^- c_h$, either works).*

*Proof.* To simplify notation, we omit writing the dependencies of the quantities on $m$. We let $\widetilde{c}_b$ the smooth extension of $c_b$ to the whole domain and subtract $\widetilde{c}_b(t^m + s\Delta t_c, x_{\underline{i}}) =: \widetilde{c}_{\underline{i}}^s$ from the scheme to obtain homogeneous Dirichlet boundary conditions and define $\widehat{c}_{\underline{i}}^s := c_{\underline{i}}^s - \widetilde{c}_{\underline{i}}^s$. $\widehat{c}_{\underline{i}}^s$ satisfies the scheme

$$D_{t_c}^+ \widehat{c}_{\underline{i}}^s - \nu_c \Delta_h \widehat{c}_{\underline{i}}^s = (\widehat{c}_{\underline{i}}^s + \widetilde{c}_{\underline{i}}^s)\boldsymbol{\Psi}_c(n_{\underline{i}}, \widehat{c}_{\underline{i}}^s + \widetilde{c}_{\underline{i}}^s) + r_c(c_{\text{supp}} - \widehat{c}_{\underline{i}}^s - \widetilde{c}_{\underline{i}}^s) - D_{t_c} \widetilde{c}_{\underline{i}}^s + \nu_c \Delta_h \widetilde{c}_{\underline{i}}^s. \quad (4.15)$$

Let us denote $H_1(\widehat{c}_{\underline{i}}^s, n_{\underline{i}}) := (\widehat{c}_{\underline{i}}^s + \widetilde{c}_{\underline{i}}^s)\boldsymbol{\Psi}_c(n_{\underline{i}}, \widehat{c}_{\underline{i}}^s + \widetilde{c}_{\underline{i}}^s) + r_c(c_{\text{supp}} - \widehat{c}_{\underline{i}}^s - \widetilde{c}_{\underline{i}}^s) - D_{t_c} \widetilde{c}_{\underline{i}}^s + \nu_c \Delta_h \widetilde{c}_{\underline{i}}^s$. By Lemmas 4.2 and 4.4 and because $\widetilde{c}_b$ is smooth, we know that $H_1(\widehat{c}_h, n_h) \in L^\infty((0,T) \times \Omega)$ uniformly for all $h > 0$. We multiply equation (4.15) by $\widehat{c}_{\underline{i}}^s$, sum over all $\underline{i}$ and rename the indices to obtain

$$\sum_{\underline{i}} \widehat{c}_{\underline{i}}^s D_{t_c}^+ \widehat{c}_{\underline{i}}^s = \sum_{\underline{i}} \left(\nu_c \widehat{c}_{\underline{i}}^s \Delta_h \widehat{c}_{\underline{i}}^s + \widehat{c}_{\underline{i}}^s H_1(\widehat{c}_{\underline{i}}^s, n_{\underline{i}})\right)$$
$$= \sum_{\underline{i}} \left(-\nu_c |\nabla_h \widehat{c}_{\underline{i}}^s|^2 + \widehat{c}_{\underline{i}}^s H_1(\widehat{c}_{\underline{i}}^s, n_{\underline{i}})\right),$$

where we have used that $\widehat{c}_{\underline{i}}^s$ satisfy homogeneous Dirichlet boundary conditions. We note that

$$\widehat{c}_{\underline{i}}^s D_{t_c}^+ \widehat{c}_{\underline{i}}^s = \frac{1}{2} D_{t_c}^+ |\widehat{c}_{\underline{i}}^s|^2 - \frac{\Delta t_c}{2} |D_{t_c}^+ \widehat{c}_{\underline{i}}^s|^2,$$

and estimate the second term,

$$\sum_{\underline{i}} |D_{t_c}^+ \widehat{c}_{\underline{i}}^s|^2 = \sum_{\underline{i}} \left|\nu_c \Delta_h \widehat{c}_{\underline{i}}^s + \widehat{c}_{\underline{i}}^s H_1(\widehat{c}_{\underline{i}}^s, n_{\underline{i}})\right|^2$$
$$\leq \frac{1}{4} \sum_{\underline{i}} |\nu_c \Delta_h \widehat{c}_{\underline{i}}^s|^2 + 16 \sum_{\underline{i}} \left|\widehat{c}_{\underline{i}}^s H_1(\widehat{c}_{\underline{i}}^s, n_{\underline{i}})\right|^2$$
$$\leq \frac{\nu_c^2}{h^2} \sum_{\underline{i}} |\nabla_h^- \widehat{c}_{\underline{i}}^s|^2 + 16 \sum_{\underline{i}} \left|\widehat{c}_{\underline{i}}^s H_1(\widehat{c}_{\underline{i}}^s, n_{\underline{i}})\right|^2$$

Since $\widehat{c}_h H_1(\widehat{c}_h, n_h) \in L^\infty((0,T) \times \Omega)$, we have

$$\frac{1}{2} D_{t_c}^+ \left(\sum_{\underline{i}} |\widehat{c}_{\underline{i}}^s|^2\right) + \nu_c \sum_{\underline{i}} |\nabla_h \widehat{c}_{\underline{i}}^s|^2$$

$$\leq \frac{\Delta t_c \nu_c^2}{2h^2} \sum_{\underline{i}} |\nabla_h \widehat{c}_{\underline{i}}^s|^2 + \sum_{\underline{i}} \widehat{c}_{\underline{i}}^s H_1(\widehat{c}_{\underline{i}}^s, n_{\underline{i}}) \left(1 + 8\Delta t_c \, \widehat{c}_{\underline{i}}^s H_1(\widehat{c}_{\underline{i}}^s, n_{\underline{i}})\right)$$

$$\leq \frac{\nu_c}{2} \sum_{\underline{i}} |\nabla_h \widehat{c}_{\underline{i}}^s|^2 + C N_x^d,$$

where $C$ is a constant not depending oh $h$ and $\Delta t_c$, by the CFL-condition (4.11). Therefore, multiplying by $\Delta t_c h^d$ and summing over all $s$, we obtain, using $N_c \Delta t_c = \Delta t$,

$$h^d \sum_{\underline{i}} \left|\widehat{c}_{\underline{i}}^{m+1}\right|^2 + h^d \Delta t_c \nu_c \sum_{s=1}^{N_c} \sum_{\underline{i}} |\nabla_h \widehat{c}_{\underline{i}}^{m,s}|^2 \leq h^d \sum_{\underline{i}} |\widehat{c}_{\underline{i}}^m|^2 + C\Delta t.$$

Using induction on $m$ in the last estimate and that $c_h$ differs from $\widehat{c}_h$ by a bounded function with bounded first differences, we obtain that $\nabla_h c_h \in L^2([0,T] \times \Omega)$ uniformly in $h > 0$. □

*Remark* 4.6. Maximum principles, nonnegativity of $q_h$ and estimates on $\nabla_h q_h$ are proved in exactly the same way.

4.2.1. *Estimates on the discrete potential $W_h$.*

**Lemma 4.7.** *We have that*

$$W_h, \nabla_h W_h, \nabla_h^2 W_h \in L^\infty([0,T]; L^2(\Omega)),$$

*uniformly in $h > 0$, where $\nabla_h := \nabla_h^+$ or $\nabla_h := \nabla_h^-$ (either works) and $\nabla_h^2 := \nabla_h^+ \nabla_h^-$ or $\nabla_h^2 := \nabla_h^- \nabla_h^+$, and*

$$W_h, \Delta_h W_h \in L^\infty((0,T) \times \Omega)),$$

*uniformly in $h > 0$ as well.*

*Proof.* The proof of this lemma can be found in [16, Lemma 4.3]. □

4.3. **Discrete entropy inequalities for $n_h$.** To prove strong convergence of the approximating sequence $\{(n_h, W_h, c_h, q_h)\}_{h>0}$, it will be useful to derive entropy inequalities for $n_h$. To this end, the following lemma will be useful:

**Lemma 4.8.** *Let $f : \mathbb{R} \to \mathbb{R}$ be a smooth convex function and assume that $\Delta t$ satisfies the CFL-condition (4.5). Denote $f_{\underline{i}}^m := f(n_{\underline{i}}^m)$ and $f_h$ a piecewise constant interpolation of it as in (4.4). Then $f_{\underline{i}}^m$ satisfies the following identity*

$$D_t f_{\underline{i}}^m = \operatorname{div}_h^- \mathbf{Q}_{\underline{i}}^m - \frac{h}{4} \sum_{j=1}^d [f''(\widehat{n}_{\underline{i}+\frac{1}{2}\mathbf{e}_j}^m)|D_j^+ W_{\underline{i}}^m||D_j^+ n_{\underline{i}}^m|^2 \qquad (4.16)$$

$$+ f''(\widehat{n}_{\underline{i}-\frac{1}{2}\mathbf{e}_j}^m)|D_j^- W_{\underline{i}}^m||D_j^- n_{\underline{i}}^m|^2] \qquad (4.17)$$

$$+ (f'(n_{\underline{i}}^m)n_{\underline{i}}^m - f_{\underline{i}}^m)\Delta_h W_{\underline{i}}^m + f'(n_{\underline{i}}^m)n_{\underline{i}}^m \mathbf{\Phi}(p_{\underline{i}}^m, c_{\underline{i}}^m, q_{\underline{i}}^m) \qquad (4.18)$$

$$+ \frac{\Delta t}{2} f''(\widetilde{n}_{\underline{i}}^{m+1/2})|D_t^+ n_{\underline{i}}^m|^2, \qquad (4.19)$$

*where*

$$\mathbf{Q}_{\underline{i}}^{m,(j)} = \frac{f_{\underline{i}}^m + f_{\underline{i}+\mathbf{e}_j}^m}{2} D_j^+ W_{\underline{i}}^m + \frac{h}{2} \frac{f'(n_{\underline{i}}^m) + f'(n_{\underline{i}+\mathbf{e}_j}^m)}{2} |D_j^+ W_{\underline{i}}^m| D_j^+ n_{\underline{i}}^m$$

$$- \frac{h^2}{4} f''(\widetilde{n}_{\underline{i}+\frac{1}{2}\mathbf{e}_j}^m) D_k^+ W_{\underline{i}}^m |D_j^+ n_{\underline{i}}^m|^2$$

*where $\widetilde{n}_{\underline{i}\pm\frac{1}{2}\mathbf{e}_j}^m, \widehat{n}_{\underline{i}\pm\frac{1}{2}\mathbf{e}_j}^m \in [\min\{n_{\underline{i}}^m, n_{\underline{i}\pm\mathbf{e}_j}^m\}, \max\{n_{\underline{i}}^m, n_{\underline{i}\pm\mathbf{e}_j}^m\}]$, and*

$$\widetilde{n}_{\underline{i}}^{m+1/2} \in [\min\{n_{\underline{i}}^m, n_{\underline{i}}^{m+1}\}, \max\{n_{\underline{i}}^m, n_{\underline{i}}^{m+1}\}]$$

*and where the term* (4.18) *is uniformly bounded and the terms* (4.16), (4.17) *and* (4.19) *satisfy for* $j = 1, \ldots, d$,

$$\frac{h^{d+1}\Delta t}{2} \sum_{m=0}^{N_t} \sum_{\underline{i}} f''(\widehat{n}^m_{\underline{i}+\frac{1}{2}\mathbf{e}_j}) |D_j^+ W_{\underline{i}}^m| |D_j^+ n_{\underline{i}}^m|^2 \leq C,$$

$$\frac{h^d \Delta t^2}{2} \sum_{m=0}^{N_t} \sum_{\underline{i}} f''(\widetilde{n}^{m+1/2}_{\underline{i}}) |D_t^+ n_{\underline{i}}^m|^2 \leq C, \tag{4.20}$$

*In particular, this implies that the piecewise constant interpolation $D_t^+ f_h$ is of the form $D_t^+ f_h = g_h + k_h$ where $g_h \in L^1([0,T] \times \Omega)$ and $k_h \in L^\infty([0,T]; W^{-1,q}(\Omega))$ for any $1 \leq q < \infty$ if $d = 2$ and for $1 \leq q \leq q^* = 2d/(d-2)$ if $d > 2$, uniformly in $h > 0$.*

*Proof.* This lemma is a slightly more general version of Lemma 4.5 in [16], but the proof is done in the exact same way. For details, check out [16, Lemma 4.5]. □

*Remark* 4.9. The preceeding lemma implies that the forward time difference of the approximation of the pressure $D_t^+ p_h = D_t^+ |n_h|^\gamma$ is of the form $D_t^+ p_h = g_h + k_h$ where $g_h \in L^1([0,T] \times \Omega)$ and $k_h \in L^\infty([0,T]; W^{-1,q}(\Omega))$ for any $1 \leq q < \infty$ if $d = 2$ and for $1 \leq q \leq q^* = 2d/(d-2)$ if $d > 2$, uniformly in $h > 0$. Using this, we have that $D_t^+ W_h = U_h + V_h$ where $U_h$ solves

$$-\mu \Delta_h U_h + U_h = g_h, \text{ and } -\mu \Delta_h V_h + V_h = k_h.$$

By Lemma B.1 in [16], we have $U_h, \nabla_h U_h \in L^1([0,T]; L^q(\Omega))$ for $1 \leq q \leq d/(d-1)$ and by standard results, $V_h, \nabla_h V_h \in L^\infty([0,T]; L^2(\Omega))$. Hence $D_t W_h, D_t \nabla_h W_h \in L^1([0,T]; L^q(\Omega)) + L^\infty([0,T]; L^2(\Omega))$.

4.4. **Passing to the limit** $h \to 0$. The estimates of the previous (sub)sections allow us to pass to the limit $h \to 0$ in a subsequence still denoted $h$ and conclude existence of limit functions

$$n_h \rightharpoonup n \geq 0, \quad \text{in } L^q([0,T] \times \Omega), 1 \leq q < \infty,$$
$$p_h \rightharpoonup \overline{p} \geq 0, \quad \text{in } L^q([0,T] \times \Omega), 1 \leq q < \infty,$$

where $p_h := n_h^\gamma$ and $0 \leq n, \overline{p} \in L^\infty([0,T] \times \Omega)$. Using the "discretized" Aubin-Lions lemma, Lemma B.1, (c.f. [16, Lemma A.1]) for $W_h, \nabla_h W_h, c_h$ and $q_h$, we obtain strong convergence of a subsequence in $L^q([0,T] \times \Omega)$ for any $q \in [0,\infty)$ in the case of $W_h$ and $1 \leq q \leq 2^*$ in the case of $\nabla_h W_h, c_h$ and $q_h$ ($2^* = 2d/(d-2)$ if $d \geq 3$ and any finite number greater than or equal to one if $d = 2$), to limit functions $W, \nabla W, c, q \in L^2([0,T]; H^1(\Omega))$. Moreover, from the estimates in Lemma 4.7 we obtain that $W \in L^\infty([0,T] \times \Omega) \cap L^\infty([0,T]; H^2(\Omega))$. Hence we have that $(n, W, \overline{p}, c, q)$ satisfy (3.4) for any test functions $\varphi, \psi_j \in C^1([0,T] \times \Omega)$, $j = 1, 2, 3$. To conclude that the limit $(n, W, \overline{p}, c, q)$ is a weak solution of (1.9), we proceed as in the previous Section 3 and show that $n_h$ in fact converges strongly. To do so, we use the discrete entropy inequality from Lemma 4.8 and Lemma 3.10, to obtain inequalities (3.9) and (3.8). Subtracting them from each other, we get (3.10). The terms on the right hand side of the resulting expression can then be estimated in the same way as in the continuous setting in Section 3.3, after noting that a discrete version of Lemma 3.11 holds (cf. [16, Lemma 4.8]). We conclude that the approximations of the cell density $n_h$ converge strongly to $n$ and that the limit is in fact a weak solution of (1.9).

## 5. Numerical examples

In this section, we test the scheme from Section 4 on several different initial data and parameters.

5.1. **Necrotic core.** We consider for the cell density $n$ the initial data

$$n_0(x) = \frac{1}{2} \exp\left(-10\left((x_1 - 0.7)^2 + x_2^2\right)\right) \\ + \frac{1}{2} \exp\left(-20\left((x_1 + 0.6)^2 + (x_2 - 0.2)^2\right)\right) \quad (5.1)$$

with homogeneous Neumann boundary conditions for both $n$ and $W$ on the domain $\Omega = [-3,3]^2$ and $h = 1/64$ with pressure law $p = n^{10}$ and $\boldsymbol{G}(p) = 1 - p$ and $\mu = 1$. That means that at time $T = 0$ two small colonies of tumor cells are present at $x = (0.7, 0)$ and $x = (-0.6, 0.2)$. We assume that there is no drug present and that the nutrient is evenly distributed over the whole domain at a level $c_0 \equiv 1$. We set $c_b \equiv c_{\text{supp}} \equiv 1$. We use $g_1(c,q) = 8(c - c_{\text{crit}})\mathbb{1}_{c > c_{\text{crit}}}$, $g_2(c,q) = 8(c_{\text{crit}} - c)\mathbb{1}_{c < c_{\text{crit}}}$ and $\boldsymbol{\Psi}_c(n,c) = -20nc$. We set the critical nutrient concentration $c_{\text{crit}} = 0.25$. The simulation results are displayed in Figures 1 (cell density) and 2 (nutrient concentration). They establish that after a rapid initial growth of the tumor cell colony due to the nutrient supply, the cell density starts decreasing again in the core of the tumor after time $t = 4$. Figures 2 exhibit that at the same time the nutrient falls below the critical level $c_{\text{crit}}$ inducing cell death in the core.

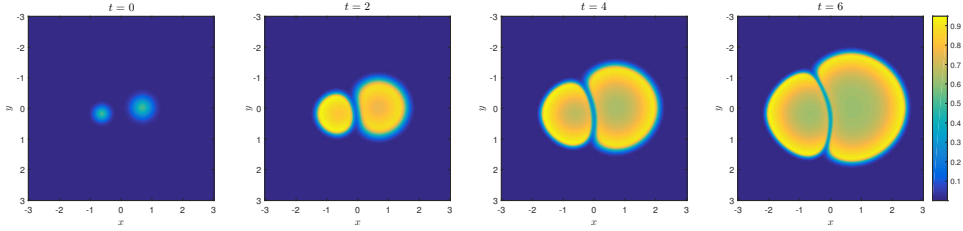

FIGURE 1. The approximations of the cell density $n$ for initial data (5.1) on $\Omega = [-3,3]^2$ with mesh width $h = 1/64$.

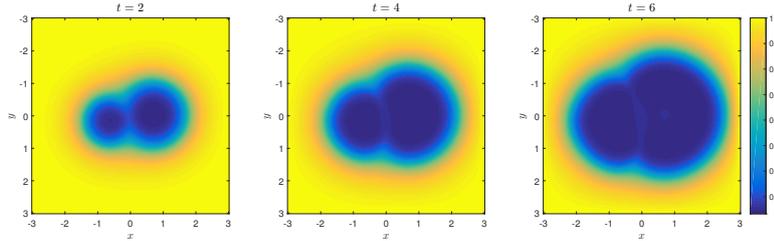

FIGURE 2. The approximations of the nutrient $c$ for constant initial data $c_0 \equiv 1$, initial data (5.1) for $n$, on $\Omega = [-3,3]^2$ with mesh width $h = 1/64$.

5.2. **With drug application.** Now we consider the same example but with the addition of a drug $q$. We set an initial drug concentration $q_0(x) \equiv 1$ constant over the domain and assume that the supply is constant over time ($q_b(t) \equiv q_{\text{supp}} \equiv 1$). We use $\boldsymbol{\Psi}_q(n,q) = -15nq$ in the diffusion equation for the drug $q$. For $g_2$, we use $g_2(c,q) = 8(c_{\text{crit}} - c)\mathbb{1}_{c < c_{\text{crit}}} + 4(q - q_{\text{crit}})\mathbb{1}_{q > q_{\text{crit}}}$ for $q_{\text{crit}} = 0$, the other functions remain as in the previous examples. The results are displayed in Figures 3, 4, 5. The numerical experiments establish that due to the drug application the initial growth of the tumor is not as rapid as in the previous example. Moreover, the cell density does not decrease as rapidly at the center of the tumor, as a result the size of the necrotic core is smaller in this case. We observe that, depending on the balance between nutrient and drug application, both shrinkage and growth of tumors occur.

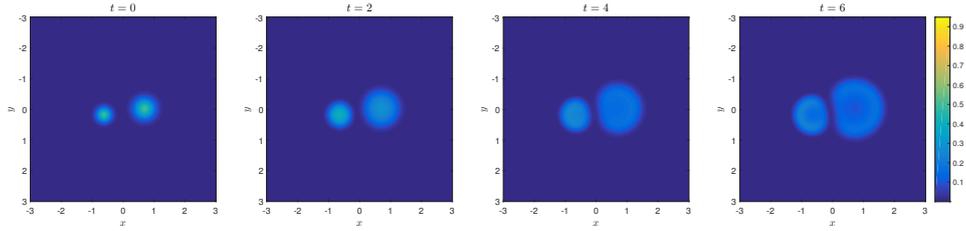

FIGURE 3. The approximations of the cell density $n$ for initial data (5.1) on $\Omega = [-3, 3]^2$ with mesh width $h = 1/64$ for the example with drug application.

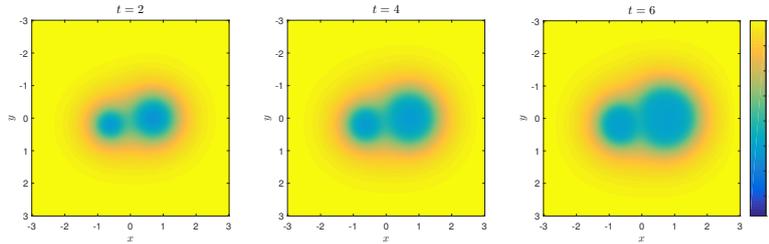

FIGURE 4. The approximations of the nutrient $c$ for constant initial data $c_0 \equiv$ and $q_0 \equiv 1$, initial data (5.1) for $n$ on $\Omega = [-3, 3]^2$ with mesh width $h = 1/64$ for the example with drug application.

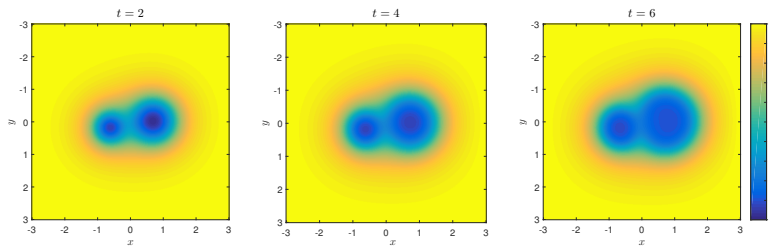

FIGURE 5. The approximations of the drug $q$ for constant initial data $c_0 \equiv$ and $q_0 \equiv 1$, initial data (5.1) for $n$ on $\Omega = [-3, 3]^2$ with mesh width $h = 1/64$.

5.3. **Shape irregularities.** In this example, we investigate the influence of the shape of the boundary of the domain on the growth of the tumor. We use the symmetric initial data

$$n_0(x) = \frac{1}{2} \exp\left(-10\left(x_1^2 + x_2^2\right)\right) \tag{5.2}$$

with homogeneous Neumann boundary conditions for both $n$ and $W$ on the domain $\Omega = [-5,5]^2$ and $h = 1/64$ with pressure law $p = n^{30}$ and $\boldsymbol{G}(p) = 1-p$ and $\mu = 0.1$. That means that at time $T = 0$ there is one small colony of tumor cells present at $x = (0,0)$. We assume that there is no drug present and that the nutrient is initially evenly distributed over the whole domain at a level $c_0 \equiv 1$. We set $c_b \equiv c_{\text{supp}} \equiv 1$. We use $g_1(c,q) = 200(c-c_{\text{crit}})\mathbb{1}_{c>c_{\text{crit}}}$, $g_2(c,q) = 200(c_{\text{crit}}-c)\mathbb{1}_{c<c_{\text{crit}}}$ and $\boldsymbol{\Psi}_c(n,c) = -20nc$. We set the critical nutrient concentration $c_{\text{crit}} = 0.5$, $r_c = 0.0001$ and $\nu_c = 5$. This means that most of the nutrient is diffused into the tumor region from the boundary of $\Omega$ and only few blood vessels are present within the tumor region (the nutrient supply $r_c c_{\text{supp}}$ is small).

The simulation results are displayed in Figures 6 (cell density) and 7 (nutrient concentration). We notice that, as time evolves, the tumor looses its spherically symmetric surface and develops irregularities in shape. It appears that the shape of the domain $\Omega$ and thereby the different levels of nutrient concentration and gradients influence this behavior. Related observations were made in [1].

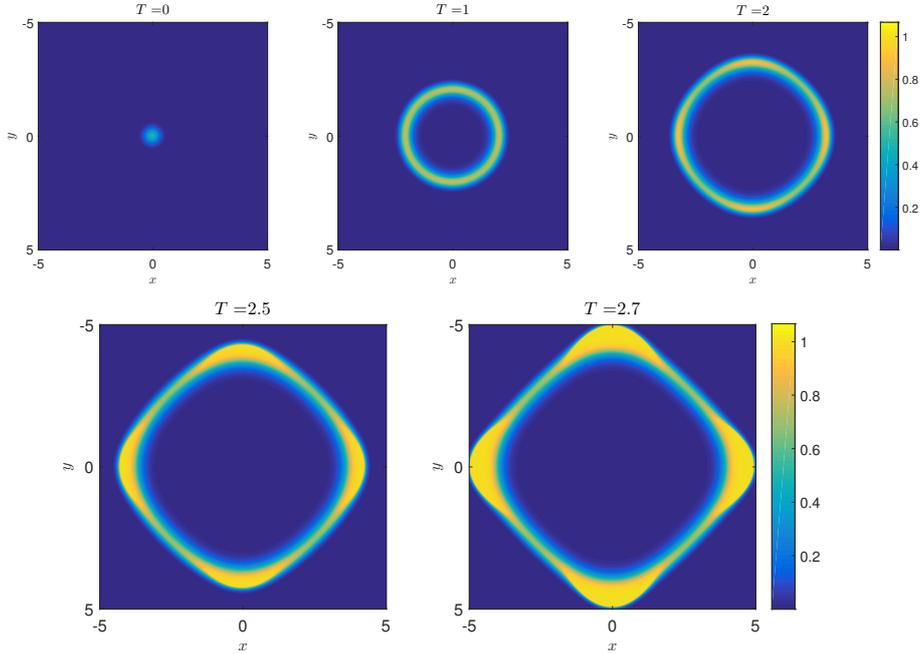

FIGURE 6. The approximations of the cell density $n$ for initial data (5.2) and parameters as in Section 5.3 on $\Omega = [-5,5]^2$ with mesh width $h = 1/64$.

5.4. **Inhomogeneous boundary conditions.** In the last experiment, we observed that the shape of the domain and so the nutrient supply from the boundary has an effect on the shape the tumor develops. To investigate this phenomenon in more detail, we will in the following experiment use inhomogeneous boundary data for the nutrient, specifically, we will use the initial data

$$c_0(x, y) = 0.8 + 0.5 \sin(0.2\,\pi\,y), \qquad (5.3)$$

for the nutrient and its restriction to the boundary $\partial\Omega$ as boundary data $c_b$. Apart from that, we use the same parameters as in Section 5.3. The results are displayed in Figures 8 (cell density) and 9 (nutrient).

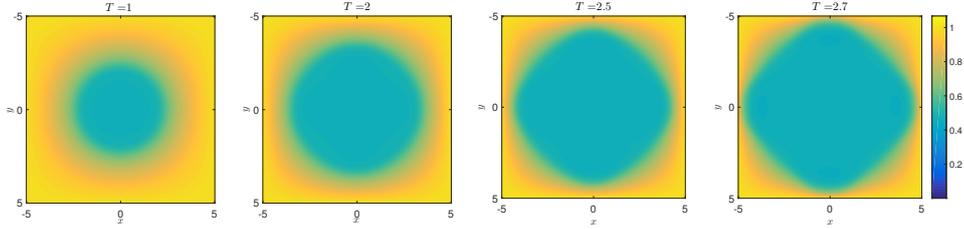

FIGURE 7. The approximations of the nutrient $c$ for initial data (5.2) and parameters as in Section 5.3 on $\Omega = [-5, 5]^2$ with mesh width $h = 1/64$ for the example with drug application.

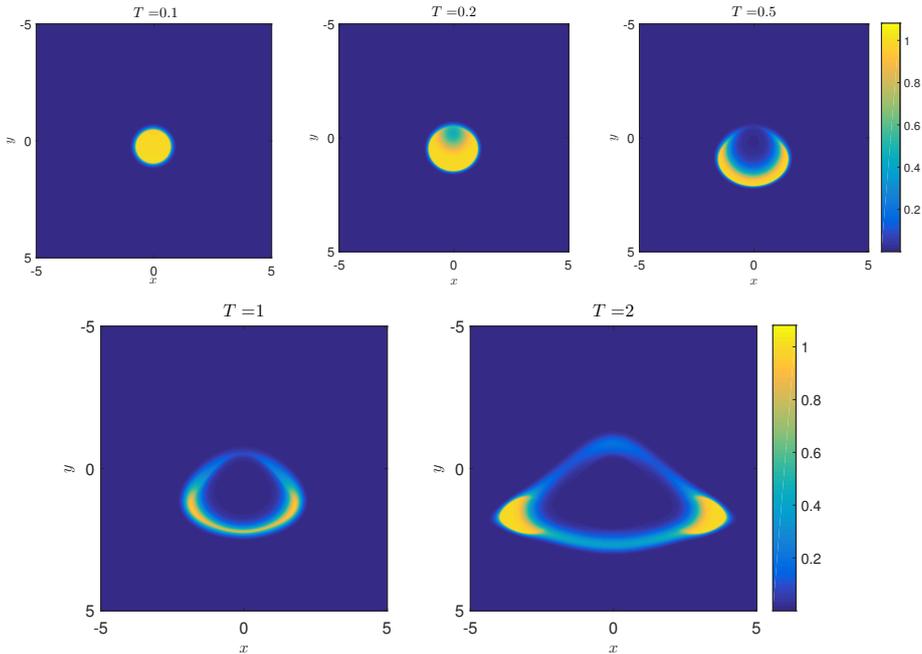

FIGURE 8. The approximations of the cell density $n$ for initial data (5.2), (5.3) and parameters as in Section 5.4 on $\Omega = [-5, 5]^2$ with mesh width $h = 1/64$.

Indeed, the inhomogeneous supply of nutrient affects the developing shape the tumor significantly. Figure 9 present the approximation of the nutrient in the case of initial data (5.2) and (5.3) respectively. The example shows that in the case of inhomogeneous boundary condition the model exhibits faster decay of the concentration of the nutrient in the center of the tumor.

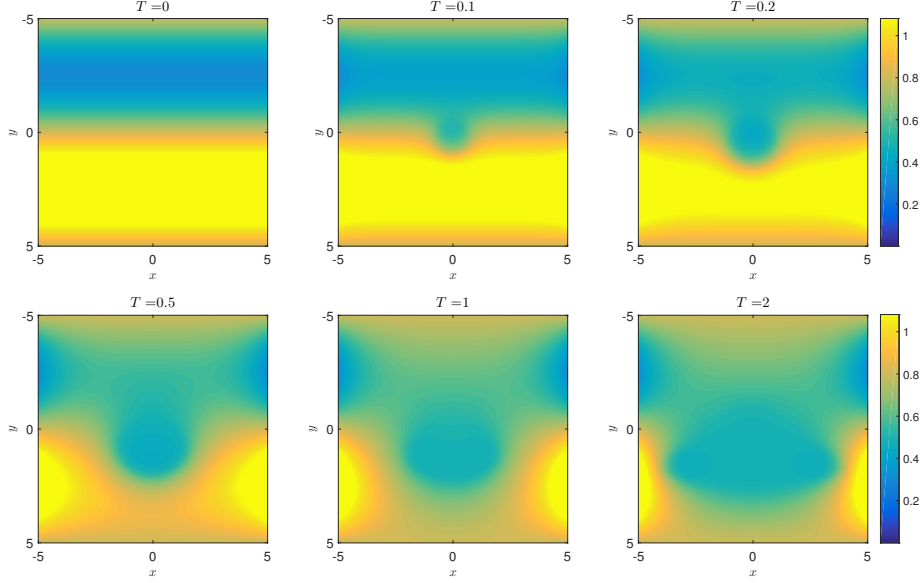

FIGURE 9. The approximations of the nutrient $c$ for initial data (5.2), (5.3) and parameters as in Section 5.4 on $\Omega = [-5, 5]^2$ with mesh width $h = 1/64$.

We conclude that inhomogeneous supply of nutrient in the tumor regions and the different level of vascularization in the tumor regime can result in shape irregularities which may lead to tumor break-up.

## APPENDIX A. PROOF OF LEMMA 3.10

*Proof.* We let $0 \leq \psi \in C_0^\infty(\mathbb{R}^{d+1})$ be a smooth, radially symmetric mollifier, i.e. $\psi(x) = \psi(-x)$ and $\int_{\mathbb{R}^{d+1}} \psi(x)dx$, with $\mathrm{supp}(\psi) \subset B_1(0)$ and denote for $\delta > 0$, $\psi_\delta(x) := \delta^{-(d+1)}\psi(x/\delta)$. Then we choose as a test function in (3.5) $\psi_\delta(s,y)\varphi(t+s, x+y)$, with $\varphi$ is compactly supported in $(\delta, T-\delta) \times \Omega^\delta$ where $\Omega^\delta$ includes all the points $x$ in $\Omega$ which have distance $d(x, \partial\Omega) > \delta$ and do a change of variables:

$$\int_0^T \int_\Omega n(t-s, x-y)\psi_\delta(s,y)\partial_t\varphi(t,x) - n(t-s, x-y)\boldsymbol{u}(t,x)\psi_\delta(s,y) \cdot \nabla\varphi(t,x)\, dxdt$$
$$= -\int_0^T \int_\Omega f(t-s, x-y)\psi_\delta(s,y)\varphi(t,x)\, dxdt.$$

Integrating in $(s, y)$, this becomes

$$\int_0^T \int_\Omega (n * \psi_\delta)(t,x)\partial_t\varphi(t,x) - (n\boldsymbol{u}) * \psi_\delta(t,x) \cdot \nabla\varphi(t,x)\, dxdt$$
$$= -\int_0^T \int_\Omega (f * \psi_\delta)(t,x)\varphi(t,x)\, dxdt.$$

We define $n_\delta := n * \psi_\delta$ and $f_\delta := f * \psi_\delta$ and choose as a test function $\varphi := b'(n_\delta)\phi$ for a smooth $\phi$ compactly supported in $(\delta, T - \delta) \times \Omega^\delta$ (which is possible since $n_\delta$ is smooth and bounded thanks to the convolution.). Then we can rewrite the last identity using chain rule as

$$\int_0^T \int_\Omega b(n_\delta)\partial_t\phi - b(n_\delta)\boldsymbol{u} \cdot \nabla\phi\, dxdt$$

$$= -\int_0^T \int_\Omega (b'(n_\delta)f_\delta + [b'(n_\delta)n_\delta - b(n_\delta)] \operatorname{div} \boldsymbol{u} + b'(n_\delta)r_\delta) \phi \, dxdt.$$

where $r_\delta := \operatorname{div}((n\boldsymbol{u}) * \psi_\delta) - \operatorname{div}(n_\delta \boldsymbol{u})$. By [11, Lemma 2.3], we have that $r_\delta \to 0$ in $L^2_{\mathrm{loc}}((0,T) \times \Omega)$ and thanks to the properties of the convolution that $b(n_\delta) \to b(n)$ almost everywhere as well as $f_\delta \to f$ a.e. when $\delta \to 0$. Thus we obtain that in the limit $\delta \to 0$, $n$ satisfies

$$\int_0^T \int_\Omega b(n)\partial_t \phi - b(n)\boldsymbol{u} \cdot \nabla \phi \, dxdt = -\int_0^T \int_\Omega (b'(n)f + [b'(n)n - b(n)] \operatorname{div} \boldsymbol{u}) \phi \, dxdt.$$

which is exactly (3.6) in the sense of distributions. $\square$

## Appendix B. Discretized Aubin-Lions Lemma

**Lemma B.1.** *Let $u_h : [0,T) \times \Omega \to \mathbb{R}^k$ be a piecewise constant function defined on a grid on $[0,T) \times \Omega$, $\Omega$ a bounded rectangular domain, satisfying*

$$\int_0^T \int_\Omega |u_h|^q + |\nabla_h u_h|^q \, dxdt \leq C \tag{B.1}$$

*for some $\infty > q > 1$, uniformly with respect to $h > 0$ and*

$$D_t u_h = A_h f_h + g_h + k_h, \tag{B.2}$$

*where $A_h$ is a first order linear finite difference operator, and $f_h, g_h, k_h : \Omega \to \mathbb{R}^{d \times k}$ are piecewise constant functions, satisfying uniformly in $h > 0$,*

$$\int_0^T \int_\Omega |f_h|^{r_1} + |g_h|^{r_2} + |k_h| \, dxdt \leq C, \tag{B.3}$$

*for some $\infty > r_1, r_2 > 1$. Then $u_h \to u$ in $L^q([0,T) \times \Omega)$.*

*Proof.* The proof can be found in [16, Lemma A.1]. $\square$

ON A TUMOR GROWTH MODEL 25[11] P.-L. Lions. *Mathematical topics in fluid mechanics. Vol. 1*, volume 3 of *Oxford Lecture Series in Mathematics and its Applications*. The Clarendon Press, Oxford University Press, New York, 1996. Incompressible models, Oxford Science Publications.
[12] A. Lunardi. *Analytic semigroups and optimal regularity in parabolic problems*. Progress in Nonlinear Differential Equations and their Applications, 16. Birkhäuser Verlag, Basel, 1995.
[13] A. Novotný and I. Straškraba. *Introduction to the mathematical theory of compressible flow*, volume 27 of *Oxford Lecture Series in Mathematics and its Applications*. Oxford University Press, Oxford, 2004.
[14] B. Perthame, F. Quirós, and J. L. Vázquez. The Hele-Shaw asymptotics for mechanical models of tumor growth. *Arch. Ration. Mech. Anal.*, 212(1):93–127, 2014.
[15] B. Perthame, M. Tang, and N. Vauchelet. Traveling wave solution of the Hele-Shaw model of tumor growth with nutrient. *Math. Models Methods Appl. Sci.*, 24(13):2601–2626, 2014.
[16] K. Trivisa and F. Weber. A convergent explicit finite difference scheme for a mechanical model for tumor growth. *arXiv preprint arXiv:1504.05982*, 2015.
[17] J. Ward and J. King. Mathematical modelling of avascular-tumour growth ii: modelling growth saturation. *Mathematical Medicine and Biology*, 16(2):171–211, 1999.
[18] J. P. Ward and J. King. Mathematical modelling of avascular-tumour growth. *Mathematical Medicine and Biology*, 14(1):39–69, 1997.
(K. Trivisa)

Department of Mathematics, University of Maryland, College Park, MD 20742, USA.
*E-mail address*: trivisa@math.umd.edu
*URL*: math.umd.edu/~trivisa

(F. Weber)

Seminar for Applied Mathematics, ETH Zürich, 8092 Zürich, Switzerland.
*E-mail address*: franziska.weber@sam.math.ethz.ch